\documentclass{amsart}

\usepackage{amsfonts}
\usepackage{amssymb}
\usepackage{mathabx}
\usepackage{amscd}
\usepackage{pictexwd,dcpic}
\usepackage{graphicx}
\usepackage{hyperref}
\hypersetup{
    colorlinks=true,
    linkcolor=blue,
    filecolor=magenta,
    urlcolor=cyan,
    citecolor=green,}
\usepackage[nameinlink,capitalize]{cleveref}

\title{Residual intersections and modules with Cohen-Macaulay Rees algebra}
\author{Alessandra Costantini}
\address{Department of Mathematics, University of California, Riverside, Riverside CA 92521, United States}
\email{alessanc@ucr.edu}
\date{}	

\newtheorem{thm}{Theorem}[section]
\newtheorem{defn}[thm]{Definition}
\newtheorem{tdefn}[thm]{Theorem and Definition}
\newtheorem{prop}[thm]{Proposition}
\newtheorem{set}[thm]{Setting}
\newtheorem{notat}[thm]{Notation}
\newtheorem{lemma}[thm]{Lemma}
\newtheorem{cor}[thm]{Corollary}

\begin{document}

\begin{abstract}
   In this paper, we consider a finite, torsion-free module $E$ over a Gorenstein local ring. We provide sufficient conditions for $E$ to be of linear type and for the Rees algebra $\mathcal{R}(E)$ of $E$ to be Cohen-Macaulay. Our results are obtained by constructing a generic Bourbaki $I$ ideal of $E$ and exploiting properties of the residual intersections of $I$.
\end{abstract}

\maketitle

\section{Introduction}

 Let $R$ be a Noetherian ring, and $E$ a finitely generated $R$-module with a rank. The Rees algebra $\,\mathcal{R}(E)\,$ of $E$ is defined as the symmetric algebra $\,\mathcal{S}(E)\,$ modulo its $R$-torsion submodule. The main goal of this work is to provide sufficient conditions for $\,\mathcal{R}(E)\,$ to be Cohen-Macaulay. 
  
 Rees algebras arise in Algebraic Geometry as homogeneous coordinate rings of blow-ups of schemes along a subscheme or as graphs of rational maps between varieties, and have been studied in connection with resolution of singularities. In many geometric situations, one is interested in Rees algebras of modules which are not ideals. For instance, the homogeneous coordinate ring of a sequence of successive blow-ups of a scheme along two or more distinct subschemes is the Rees algebra of a direct sum of two or more ideals. Moreover, given a subvariety $X$ of an affine space, the conormal variety of $X$ and the graph of the Gauss map from $X$ are projective spectra of the Rees algebras of modules which are not ideals. 

 When $E=I$ is an $R$-ideal, the Rees algebra $\mathcal{R}(I)$ coincides with the subalgebra 
    $$ R[It] = \oplus_{j \geq 0\,} I^j t^j $$
 of the polynomial ring $R[t]$. If $R$ is Cohen-Macaulay and $I$ has positive height, the first step in order to study the Cohen-Macaulay property of $\mathcal{R}(I)$ is to investigate the Cohen-Macaulayness of the associated graded ring 
    $$ \mathcal{G}(I) \coloneq \oplus_{j \geq 0\,} I^{j}/I^{j+1} \cong \mathcal{R}(I) /I\mathcal{R}(I).$$
 In fact, in this case $\,\mathcal{G}(I)$ is Cohen-Macaulay whenever $\,\mathcal{R}(I)$ is \cite{HuGrI}, and although the converse is not true in general, it holds if furthermore some numerical conditions are satisfied \cite{{IkedaTrung},{JK},{SUVDegPolyRel},{UVeqLinPres},{UlrichExpRedNo},{PU99}}. In turn, the Cohen-Macaulay property of $\,\mathcal{G}(I)$ is well-understood (see for instance \cite{{VV},{Sally},{RVpaper},{Rossi},{HSV},{HSVII},{HuHu93},{JU},{GNN},{PXie}}).

 By analogy with the case of ideals, one defines the powers of an $R$-module $E$ as $\,E^j \coloneq [\mathcal{R}(E)]_j, \,$ where $\,[\mathcal{R}(E)]_j$ denotes the degree-$j$ component of $\mathcal{R}(E)$. However, since $E^{n+1}$ is not contained in $E^n$, there is no module analogue for the associated graded ring. Hence, the study of the Cohen-Macaulay property of Rees algebras of modules requires a completely different approach than in the case of ideals and much less is known about it (see for instance \cite{{HSVlintype},{SUVjacduals},{SUV2003},{SUV2012},{Lin}}).
  
 A systematic approach to the problem was provided by Simis, Ulrich and Vasconcelos in 2003 \cite{SUV2003}. Given a finite $R$-module $E$ with a rank, they construct a larger ring $S$ which is a generic extension of $R$ and an ideal $I$ of $S$ so that $\,\mathcal{R}(E)\,$ is Cohen-Macaulay if and only if $\,\mathcal{R}(I)\,$ is. $I$ is called a \emph{generic Bourbaki ideal} of $E$ (see \cref{PrelimSection}). The problem of determining a sufficient condition for $\, \mathcal{R}(E)\,$ to be Cohen-Macaulay is then ultimately reduced to that of crafting assumptions on $E$ so that $I$ has a Cohen-Macaulay Rees algebra, and one could virtually attempt to extend any result known on the Cohen-Macaulayness of Rees algebras of ideals to the case of modules. However, this idea is often not easy to put in practice, since very often transferring properties from $E$ to $I$ and backwards is only possible if $\,\mathcal{R}(E)\,$ yields a deformation of $\, \mathcal{R}(I)$, which is not always guaranteed.

 The main goal of this work is to provide a sufficient condition for the Cohen-Macaulayness of the Rees algebra $\,\mathcal{R}(E)$ of a module $E$, under the assumption that finitely many powers of $E$ have sufficiently high depth. More specifically, we aim to generalize the following theorem, proved independently by Johnson and Ulrich \cite[3.1 and 3.4]{JU} and Goto, Nakamura and Nishida \cite[1.1 and 6.3]{GNN} for the Rees algebra of ideals. The assumptions that $I$ satisfies $\,G_{\ell}$ and the Artin-Nagata property $\,AN_{\ell-k-1}$ indicate an overall nice behavior of the residual intersections of $I$; we refer the reader to \cref{PrelimSection} for detailed definitions.

  \begin{thm} [\cite{GNN, JU}] \label{JUGNN} \hypertarget{JUGNN}{}
    Let $R$ be a local Cohen-Macaulay ring of dimension $d$ with infinite residue field. Let $I$ be an $R$-ideal with analytic spread $\ell$, height $g$ and reduction number $r$, and let $k$ be an integer. Assume that $g\geq 1$, $1 \leq k \leq \ell-1$, $r \leq k$ and $I$ satisfies $\,G_{\ell}$ and $\,AN_{\ell-k-1}$. If $\,\mathrm{depth}_{\,}I^j \geq d-\ell+k-j+1\,$ for $\,1\leq j \leq k$, then $\,\mathcal{R}(I)$ is Cohen-Macaulay.
  \end{thm} 
  
 From a technical point of view, the most challenging aspect of the problem is that it is not clear whether a module analogue for the Artin-Nagata property would be well-defined or behave well under passage to a generic Bourbaki ideal. Moreover, it is not guaranteed that similar assumptions on the depths of powers of $E$ would transfer to the powers of a generic Bourbaki ideal $I$ of $E$. Nevertheless, when $R$ is Gorenstein these technical obstructions can be overcome through a deep investigation of the linear type property of $E$ and of the residual intersections of $I$.

 We now describe the contents of this paper.
 
 In \cref{PrelimSection} we review the main properties of Rees algebras of modules and of generic Bourbaki ideals. We also recall basic definitions and results on residual intersections of ideals that will be used throughout the paper. In particular, \cref{LinTypeExt}, which follows from work of Chardin, Eisenbud and Ulrich \cite{CEU}, will be crucial in order to prove our main results.
 
 In \cref{SectionLinType} we discuss the linear type property of a module. Recall that a module $E$ is said to be of \emph{linear type} if its Rees algebra $\,\mathcal{R}(E)$ can be expressed as a quotient of a polynomial ring $\,R[T_1, \ldots, T_n]\,$ modulo an ideal $\mathcal{J}$ generated by linear forms in $\,R[T_1, \ldots, T_n]$. Our first main result, \cref{myLinType}, shows that (under mild assumptions) $E$ is of linear type if $\,\mathrm{Ext}_{R}^i(E^j, R)=0\,$ for finitely many given values of $i$ and $j$. This is obtained by showing that a generic Bourbaki ideal $I$ of $E$ has good residual intersection properties as in \cref{LinTypeExt}. We also show that the assumptions of \cref{myLinType} are satisfied in a wide range of situations, and provide a general technique to construct modules of linear type (see \cref{LinTypeSub}).

 \cref{SectionCMness} is dedicated to the Cohen-Macaulay property of Rees algebras of modules. Our second main result, \cref{GenGNN}, shows that (under mild assumptions) the Rees algebra of $E$ is Cohen-Macaulay if the depths of finitely many powers of $E$ are sufficiently high, recovering \cref{JUGNN} in the case when $E=I$ is an ideal. The theorem is proved by showing that $E$ is of linear type up to a certain codimension using \cref{myLinType}, and allows us to identify several examples of modules with Cohen-Macaulay Rees algebra. A similar approach to the problem had previously appeared in work of Lin \cite{Lin}, for more restrictive classes of modules. Our methods allow us to recover all of Lin's main results, with proofs that are logically independent of her work. In particular, Lin's \cite[3.4]{Lin} is a special case of \cref{GenGNN} for modules with maximal reduction number (see \cref{Lin3.4}). Moreover, \cref{thmIdealModules} for ideal modules is a minor variation of \cite[4.3]{Lin}, and  \cref{herLinType} improves another of Lin's results \cite[3.1]{Lin} regarding modules of linear type.

 Finally, when $(R, \mathfrak{m})$ is local \cref{GenGNN} and \cref{thmIdealModules} allow us to provide sufficient conditions for the \emph{special fiber ring} $\,\mathcal{F}(E) \coloneq \mathcal{R}(E) \otimes_R R/\mathfrak{m}\,$ of a module $E$ to be Cohen-Macaulay, recovering previous work of Corso, Ghezzi, Polini and Ulrich \cite{CGPU}. These results are obtained by establishing a general technique to study the Cohen-Macaulay property of the special fiber ring of a module, and will be included in a separate article \cite{myFiberCone} (see also the author's thesis \cite{myThesis}).

\subsection*{Acknowledgements}
 The content of this article was part of the author's Ph.D. thesis. The author is deeply grateful to her advisor Bernd Ulrich for his insightful comments and observations on the preliminary versions of the results here presented, and for his valuable guidance through the background literature on residual intersections of ideals.

\section{Preliminaries} \label{PrelimSection}
  
 In this section we review basic properties of Rees algebras of modules and generic Bourbaki ideals, and recall some facts about residual intersections of ideals.

\subsection{Modules with a rank and Rees algebras}
 
 Throughout this paper, $R$ will be a Noetherian ring and $E$ will be a finite $R$-module with a rank. Recall that a finite $R$-module $E$ has a \emph{rank}, $\mathrm{rank}_{\,}E =e$, if $\,E \otimes_R \mathrm{Quot}(R) \cong (\mathrm{Quot}(R))^e, \,$ or, equivalently, if $\,E_{\mathfrak{p}} \cong R_{\mathfrak{p}} ^e\,$ for all $\mathfrak{p} \in \mathrm{Ass}(R)$. Moreover, we will often assume that $E$ is \emph{orientable}. This means that $E$ has a rank $e >0$ and the $e$th exterior power of $E$ satisfies $\, (\bigwedge^{e} E)^{\ast \ast} \cong R, \,$ where $(-)^{\ast}$ denotes the functor $\mathrm{Hom}_R (-, R))$. 

 These notions are not very restrictive, since every module admitting a finite free resolution has a rank and every module of finite projective dimension is orientable. Moreover, every finite $R$-module over an integral domain has a rank, and every finite module over a UFD is orientable. Modules of rank 1 play a special role in our work. In fact, $E$ is a torsion-free module of rank 1 if and only if it is isomorphic to an $R$-ideal $I$ of positive grade.

 A module $E$ of rank $e$ is said to satisfy \emph{condition $G_s$} if $\,\mu(E_{\mathfrak{p}}) \leq \mathrm{dim}R_\mathfrak{p} -e +1\,$ for every $\mathfrak{p} \in \mathrm{Spec}(R)$ with $\,1\leq \mathrm{dim}_{\,}R_{\mathfrak{p}} \leq s-1, \,$. If the same condition holds for all $s$, then $E$ is said to satisfy $G_{\infty}$. Notice that for an ideal $I$ the $G_s$ condition can be restated as $\, \mu(I_{\mathfrak{p}}) \leq \mathrm{dim}R_\mathfrak{p}$ for every $\mathfrak{p} \in V(I)$ with $\mathrm{dim}R_{\mathfrak{p}} \leq s-1$. Moreover, if $R$ has dimension $d$, then $I$ is $G_{\infty}$ if and only if it is $G_{d+1}$.

 Now, let $\,R^s \stackrel{\varphi}\longrightarrow R^n \twoheadrightarrow E\,$ be any presentation of $E$. The \emph{symmetric algebra} of $E$ is 
   $$ \mathcal{S}(E) = R[T_1, \dots ,T_n]/ (\ell_1, \dots, \ell_s), $$
 where $\,\ell_1, \dots, \ell_s\, $ are linear forms in $R[T_1, \dots ,T_n]$ so that 
    $$ [T_1, \dots ,T_n] \cdot \varphi= [\ell_1, \dots, \ell_s]\footnote{\,This definition is independent of the choice of $\varphi$. In fact, $\mathcal{S}(E)$ can be equivalently described through a universal property, as a quotient of the tensor algebra of $E$ modulo the commutators of elements of degree 1 (see \cite[Section 1.6]{BH}).}. $$
     
 If $E$ has a rank, the \emph{Rees algebra} $\mathcal{R}(E)$ of $E$ is defined as the quotient of $\mathcal{S}(E)$ modulo its $R$-torsion submodule. In particular, 
   $$ \mathcal{R}(E) = R[T_1, \dots ,T_n]/ \mathcal{J}\,$$
 for some ideal $\mathcal{J}$, called the \emph{defining ideal} of $\mathcal{R}(E)$. The module $E$ is said to be of \emph{linear type} if $\,\mathcal{S}(E)$ is torsion-free, since in this case $\,\mathcal{J}=(\ell_1, \ldots, \ell_s)\,$ is an ideal of linear forms. Modules of linear type will be the main focus of \cref{SectionLinType}.
  
 When $E=I$ is an $R$-ideal, $\,\mathcal{R}(I)$ coincides with the subalgebra $\,R[It]= \oplus_{j \geq 0\,} I^j t^j\,$ of the polynomial ring $R[t]$. By analogy with the case of ideals, the powers of an $R$-module $E$ are then defined as 
  $$E^j \coloneq [\mathcal{R}(E)]_j.$$
  
 A more general definition of Rees algebra which also applies to modules not having a rank was given by Eisenbud, Huneke and Ulrich in \cite{EHU}. In the same paper, they also defined the notions of minimal reduction and analytic spread of a module $E$ over a Noetherian local ring (see \cite[2.1 and 2.3]{EHU}). A \emph{minimal reduction} of $E$ is a minimal submodule $\,U \subseteq E\,$ so that $\,E^{r+1} = U E^r\,$ for some integer $\,r \geq 0$. The least such integer $r$ is denoted by $r_U(E)$, and the \emph{reduction number} of $E$ is 
           $ \, \displaystyle{ r(E) \coloneq \mathrm{min} \, \{r _U(E) \, | \, U \mathrm{\,is \; a \; minimal \; reduction \; of\, } E\}.}$
 Finally, the \emph{analytic spread} $\ell(E)$ of $E$ is the Krull dimension of the \emph{special fiber ring} $\, \mathcal{F}(E) \coloneq \mathcal{R}(E) \otimes_R k$, where $k$ is the residue field of $R$. The following result extends the main properties of the analytic spread of an ideal to modules (see \cite[2.3]{EHU} and \cite[2.2 and 2.3]{SUV2003}).
 
  \begin{prop} \label{anspread} \hypertarget{anspread}{}
   Let $R$ be a Noetherian local ring with $\mathrm{dim}R=d$ and residue field $k$. Let $E$ be a finite $R$-module. 
     \begin{itemize}
       \item[$(a)$] If $k$ is infinite, all minimal reductions of $E$ have the same minimal number of generators, equal to $\ell(E)$. Moreover, $\,\ell(E) \leq \mu(E)$, and equality holds if and only if $E$ has no proper reductions.
       \item[$(b)$] If $k$ is infinite, any submodule $U$ generated by $\ell(E)$ general elements in $E$ is a minimal reduction of $E, \,$ with $\,r_U(E) = r(E)$. 
       \item[$(c)$] Suppose that $d>0$ and that $E$ has a rank, $\, \mathrm{rank}\,E=e$. Then 
            $$ e \leq \ell(E) \leq d +e -1 = \mathrm{dim}\, \mathcal{R}(E) -1. \,$$ 
    \end{itemize}
 \end{prop}

\subsection{Generic Bourbaki ideals}
 
 We now recall the definition and main properties of generic Bourbaki ideals, introduce by Simis, Ulrich and Vasconcelos in \cite{SUV2003}. 

 \begin{notat} \label{NotationBourbaki} \hypertarget{NotationBourbaki}{}
 $($\cite[3.3]{SUV2003}$)$.
   \emph{Let $R$ be a Noetherian ring, $E=Ra_1 + \dots + Ra_n$ a finite $R$-module with $\mathrm{rank}_{\,}E=e>0$. Let $\displaystyle Z= \{ Z_{ij} \, | \, 1 \leq i \leq n$, $1\leq j \leq e-1 \}$ be a set of indeterminates, and denote 
      $$ R'\coloneq R[Z], \quad E'\coloneq E \otimes_R R', \quad x_j \coloneq \sum_{i=1}^n Z_{ij} a_i \in E'\quad \mathrm{and} \quad F'\coloneq \sum_{j=1}^{e-1} R' x_j. $$
  If $R$ is local with maximal ideal $\mathfrak{m}$, let $\,\displaystyle R'' \coloneq R(Z)= R[Z]_{\mathfrak{m}R[Z]}$ and similarly denote $\,\displaystyle E''\coloneq E \otimes_R R''\,$ and $\, \displaystyle F''\coloneq F' \otimes_{R'} R''$.}
 \end{notat}

 \begin{thm} \label{existBourbaki} \hypertarget{existBourbaki}{}
 $($\cite[3.1, 3.2 and 3.4]{SUV2003}$)$.
   Let $R$ be a Noetherian ring, $E$ a finite $R$-module with $\mathrm{rank}_{\,}E=e>0$. Assume that $E$ is torsion-free and that $E_{\mathfrak{p}}$ is free for all $\mathfrak{p} \in \mathrm{Spec}(R)$ with $\mathrm{depth}_{\,}R_\mathfrak{p} \leq 1$. 
   \begin{itemize}
      \item[$(a)$] For $R'$, $E'$ and $F'$ as in \cref{NotationBourbaki}, $\,F'$ is a free $R'$-module of rank $e-1$ and $E'/F'$ is isomorphic to an $R'$-ideal $J$ with $\mathrm{grade}_{\,}J >0$. Moreover, $J$ can be chosen to have grade at least $\,2\,$ if and only if $E$ is orientable.
      \item[$(b)$] If $E$ satisfies $G_s$, then so does $J$.
      \item[$(c)$] If $\mathrm{grade}(J) \geq 3$, then $\,E \cong R^{e-1} \oplus L\,$ for some $R$-ideal $L$. In this case, $LR' \cong J$.
    \end{itemize}
 \end{thm}   
 
 \begin{tdefn} \label{tdefBourbaki} \hypertarget{tdefBourbaki}{}
 $($\cite[3.2 and 3.4]{SUV2003}$)$.
     \em{Let $R$ be a Noetherian local ring, and $E$ a finite $R$-module with $\mathrm{rank}_{\,}E=e>0$. With the assumptions of \cref{existBourbaki}, $\,E''/F''$ is isomorphic to an $R''$-ideal $I$, called a \emph{generic Bourbaki ideal} of $E$. Moreover, if $K$ is another ideal constructed this way using variables $Y$, then the images of $I$ and $K$ in $S=R(Z,Y)$ coincide up to multiplication by a unit in $\mathrm{Quot}(S)$, and are equal whenever $I$ and $K$ have grade at least 2.}
     \end{tdefn} 
     
 \begin{thm} \label{MainBourbaki} \hypertarget{MainBourbaki}{}
 $($\cite[3.5(a) and 3.8]{SUV2003}$)$.
   In the setting of \cref{NotationBourbaki} and with the assumptions of \cref{existBourbaki} and \cref{tdefBourbaki}, the following statements are true.
     \begin{itemize}
       \item[$(a)$] $\mathcal{R}(E)$ is Cohen-Macaulay if and only if $\,\mathcal{R}(I)$ is Cohen-Macaulay.
       \item[$(b)$] If $\,\mathrm{grade} \, \mathcal{R}(E)_+ \geq e$, then $\,\mathcal{R}(J) \cong \mathcal{R}(E')/(F')$. 
       \item [$(c)$] $E$ is of linear type with $\,\mathrm{grade} \, \mathcal{R}(E)_+ \geq e\,$ if and only if $I$ is of linear type, if and only if $J$ is of linear type.
    \end{itemize}
   \end{thm} 
 

 Thanks to \cref{MainBourbaki}, questions about the linear type property of $E$ or the Cohen-Macaulay property of $\,\mathcal{R}(E)$ are reduced to the case of ideals. Notice that in order for the conditions in \cref{MainBourbaki}(a) and (c) to hold one must have that $\,\mathcal{R}(E'')/(F'') \cong \mathcal{R}(I)\,$ and $\,x_1, \ldots, x_{e-1}\,$ form a regular sequence on $\,\mathcal{R}(E'')$ (see \cite[3.5(b)]{SUV2003}). This means that $\,\mathcal{R}(E'')$ is a \emph{deformation} of $\,\mathcal{R}(I)$. In fact, one has the following crucial technical result.

 \begin{thm} \label{SUV3.11} \hypertarget{SUV3.11}{}
 $($\cite[3.11]{SUV2003}$)$.
   Let $R$ be a Noetherian ring, $E$ a finite $R$-module with $\mathrm{rank}_{\,}E=e>0$. Let $\,0 \to F \to E \to I \to 0\,$ be an exact sequence where $F$ is a free $R$-module with free basis  $x_1, \ldots, x_{e-1}$ and $I$ is an $R$-ideal. The following are equivalent.
   \begin{itemize}
       \item[$(a)$] $\mathcal{R}(E)/(F)$ is $R$-torsion free;
       \item[$(b)$] $\mathcal{R}(E)/(F) \cong \mathcal{R}(I)$;
       \item[$(c)$] $\mathcal{R}(E)/(F) \cong \mathcal{R}(I)\,$ and $\,x_1, \ldots, x_{e-1}\,$ form a regular sequence on $\,\mathcal{R}(E)$.
    \end{itemize}
   Moreover, if $I$ is of linear type, then so is $E$ and the equivalent conditions above hold.  
 \end{thm}
 
 The construction described in \cref{NotationBourbaki} and \cref{tdefBourbaki} can be thought of as the $(e-1)$-st step of an iterative construction, where at each step the ring is extended by adjoining $n$ generic variables. This will be relevant for the proof of \cref{herLinType}. In particular, we will make use of the following result.

  \begin{thm} \label{SUV3.6} \hypertarget{SUV3.6}{}
   $($\cite[3.6 and 3.8]{SUV2003}$)$
  Let $R$ be a Noetherian ring, $E=Ra_1 + \dots + Ra_n$ a finite $R$-module with $\mathrm{rank}_{\,} E =e \geq 2$, and let $Z_1, \ldots Z_n$ be indeterminates. Denote
    $$\,\widetilde{R} \coloneq R[Z_1, \ldots, Z_n], \quad \widetilde{E} \coloneq E \otimes_R \widetilde{R}, \quad x \coloneq \sum_{i=1}^n Z_{i} a_i \in \widetilde{E} \quad \mathrm{and} \quad \overline{\mathcal{R}} \coloneq \mathcal{R}(\widetilde{E})/(x). \,$$ 
    Then
     \begin{itemize}
       \item[$($a$)$] $x$ is regular on $\mathcal{R}(\widetilde{E})$.
       \item[$($b$)$] The natural epimorphism $\,\displaystyle{\pi \colon \,  \overline{\mathcal{R}} \twoheadrightarrow \mathcal{R}(\widetilde{E}/\widetilde{R}x)\,}$ has kernel $\,\mathrm{ker}(\pi)= H^0_{\overline{\mathcal{R}}_{+}}(\overline{\mathcal{R}})$.
       \item[$($c$)$] If $\mathrm{grade}_{\,}\mathcal{R}(E)_{+} \geq 2$, then $\pi$ is an isomorphism.
     \end{itemize}
 \end{thm}
 
 The next proposition examines the connection between the analytic spread and reduction number of $E$ and those of a generic Bourbaki ideal $I$ of $E$.

 \begin{prop} \label{anspreadBourbaki} \hypertarget{anspreadBourbaki}{} $($\cite[3.10]{SUV2003}$)$.
   Let $R$ be a Noetherian local ring with $\mathrm{dim}R=d$, and let $E$ be a finite $R$-module with $\mathrm{rank}_{\,}E=e$. If $E$ admits a generic Bourbaki ideal $I$, then  $\,\ell(I)=\ell(E)-e+1$. Moreover, $\,r(I) \leq r(E)\,$ if the residue field of $R$ is infinite.
 \end{prop}

\subsection{Residual intersections}
 
 Most of the results proved in this paper rely on properties of the residual intersections of generic Bourbaki ideals of the modules under examination. 
 
 Let $R$ be a Cohen-Macaulay ring and let $I$ be an $R$-ideal with $\mathrm{ht}(I)=g$, and $\,s \geq g$ an integer. Recall that a proper ideal $K$ is an \emph{s-residual intersection} of $I\,$ if there exists an ideal $J \subseteq I$ such that $K= J \colon I$, $\mu(J) \leq s\,$ and $\,J_{\mathfrak{p}} = I_{\mathfrak{p}}\,$ for every prime ideal $\mathfrak{p}$ with $\,\mathrm{dim}R_{\mathfrak{p}} \leq s-1$. $K$ is called a \emph{geometric s-residual intersection} of $I$ if in addition $\,J_{\mathfrak{p}} = I_{\mathfrak{p}}\,$ for every $\mathfrak{p} \in V(I)$ with $\mathrm{dim}R_{\mathfrak{p}} =s$. 
 
 Residual intersections $K$ of $I$ always exist if $I$ satisfies $G_s$ \cite[1.5 and 1.6(a)]{Ulrich}, and can often be constructed using minimal reductions $J$ of $I$ (see \cite[3.1]{PXie} or \cite[2.7]{JU}). For applications to the study of Rees algebras it is often useful to consider also the following property.
    
 \begin{defn} [\cite{AN}, \cite{Ulrich}]
   \em{Let $R$ be a Cohen-Macaulay local ring, $I$ an $R$-ideal with $\mathrm{ht}(I)=g$, and $\,s \geq g$ an integer. $I$ is said to satisfy the \emph{Artin-Nagata property} $AN_s$ if for all $g \leq i \leq s$ and every $i$-residual intersection $K$ of $I$, $R/K$ is Cohen-Macaulay. If the same property holds for every geometric $i$-residual intersection, then $I$ is said to satisfy $AN_s^{-}$.}
 \end{defn}
 
 At least when $R$ is Gorenstein, the Artin-Nagata condition $AN_s$ is satisfied by a large class of ideals, including perfect ideals of height 2, perfect Gorenstein ideals of height 3, complete intersections, licci ideals and strongly Cohen-Macaulay ideals. Notice also that $AN_s$ is satisfied whenever $s \leq g-1$, or when $s=g$, $R$ is Gorenstein and $R/I$ is Cohen-Macaulay (see \cite{PS}). Another sufficient condition is given by the following result, due to Ulrich.
     
 \begin{thm}$($\cite[2.9]{Ulrich}$)$. \label{CreateAN} \hypertarget{CreateAN}{}
    Let $R$ be a local Gorenstein ring with $\mathrm{dim}(R)=d$, $I$ an $R$-ideal with $\mathrm{ht}(I)=g$. Assume that $I$ satisfies $G_s$ for some $s\geq g$. If $\mathrm{depth}(I^j) \geq d-g-j+2$ for $1 \leq j \leq s-g+1$, then $I$ satisfies $AN_s$.
 \end{thm}
 
 There is an interesting connection between the Artin-Nagata condition and the linear type property of $I$, as explained in the following result (see also \cite[1.8 and 1.9]{Ulrich}). The proof follows from the paper of Ulrich \cite{Ulrich}, but we write it here for lack of a specific reference.
 
 \begin{thm} \label{LinTypeDepth}\hypertarget{LinTypeDepth}{}
   Let $R$ be a Gorenstein local ring of dimension $d$. Let $I$ be an $R$-ideal with $\,\mathrm{ht}_{\,}I=g \geq 1\,$ and $\,\ell(I)=\ell$. Assume that $I$ is $G_{\ell+1}$ and that $\, \mathrm{depth}_{\,}I^j \geq d-g-j+2\,$ for $\,1 \leq j \leq \ell-g$. Then, $I$ is of linear type and $\,\mathcal{R}(I)\,$ is Cohen-Macaulay.
 \end{thm}
   
   \emph{Proof}. Without loss of generality, we may assume that the residue  field of $R$ is infinite. Let $J$ be a minimal reduction of $I$ generated by $\,\ell\,$ general elements. Since $I$ satisfies $G_{\ell+1}, \,$ then $\,\mathrm{ht}(J \colon_{\!R\,} I) \geq \ell+1$. This is clear if $J=I$. If $\,J \subsetneq I$, then from \cite[1.6]{Ulrich} it follows that $J \colon_{\!R\,} I$ is a geometric $\ell$-residual intersection of $I$. Hence, since $I$ and $J$ have the same radical, $\,\mathrm{ht}(J \colon_{\!R\,} I) \geq \ell+1$, as claimed. 
    
   Now, by \cref{CreateAN} the assumption on the depths of powers of $I$ imply that $I$ satisfies $AN_{\ell -1}$. Hence, if $\,J \colon_{\!R\,} I\,$ is a proper ideal, from \cite[1.7]{Ulrich} it follows that $\,\mathrm{ht}(J \colon_{\!R\,} I) =\ell$, a contradiction. So, it must be that $I=J$. This means that $I$ has no proper reductions, whence $\ell=\mu(I)$. In particular, $I$ satisfies sliding depth by \cite[1.8(c)]{Ulrich}). Hence, by \cite[6.1]{HSVII} it only remains to show that $I$ satisfies $G_{\infty}$.
   
   Notice that for every $\,\mathfrak{p} \in \mathrm{Spec}(R)\,$ one has $\,\ell(I_{\mathfrak{p}}) \leq \ell(I)$. By \cite[1.10(a)]{Ulrich} it then follows that $I_{\mathfrak{p}}$ satisfies $\,G_{\ell(I_{\mathfrak{p}})}$ and $\,AN_{\ell(I_{\mathfrak{p}})-1}$. Hence, we can repeat the argument above to deduce that $\,\mu(I_{\mathfrak{p}})=\ell(I_{\mathfrak{p}}) \leq \mathrm{dim}\,R_{\mathfrak{p}}\,$ for every $\,\mathfrak{p} \in \mathrm{Spec}(R)$. That is, $I$ satisfies $G_{\infty}$ and the proof is complete. $\blacksquare$ \\

 Thanks to work of Chardin, Eisenbud and Ulrich \cite{CEU}, the condition on the depths of powers of $I$ in \cref{LinTypeDepth} can be weakened to the requirement that certain Ext modules vanish, to produce ideals of linear type that do not necessarily satisfy the Artin-Nagata property of \cref{CreateAN} (see \cite[5.2]{CEU}). More precisely, one has the following result.
 
 \begin{thm} \label{LinTypeExt} \hypertarget{LinTypeExt}{}
   Let $R$ be a local Gorenstein ring, $d=\mathrm{dim}_{\,}R$, $I$ an $R$-ideal with $\mathrm{ht}(I)=g \geq 1$ and $\ell(I)=\ell$. Assume that $I$ is $G_{\ell+1}$ and that $\, \mathrm{Ext}_R^{\,g+j-1}(I^j, R)=0$ for $1 \leq j \leq \mathrm{min} \, \{\ell-g, \,d-g-1 \}$. Then, $I$ is of linear type.
 \end{thm}

   \emph{Proof}. Without loss of generality, we may assume that the residue field of $R$ is infinite. Let $J$ be a minimal reduction of $I$ generated by $\ell$ general elements. Since $I$ satisfies $G_{l+1}$, as in the proof of \cref{LinTypeDepth} it follows that $\,\mathrm{ht}(J \colon I) \geq \ell+1$.
 
   Now, if $\,J \colon_{\!R\,} I\,$ is a proper ideal, one must have $\ell \leq d-1$. Hence, the assumption that $\, \mathrm{Ext}_R^{\,g+j-1}(I^j, R)=0\,$ for $\,1 \leq j \leq d-g-1\,$ implies that $\,\mathrm{ht}(J \colon_{\!R\,} I) =\ell\,$ by \cite[4.1 and 3.4]{CEU}. This is a contradiction, so it must be that $I=J$, whence $\ell=\mu(I)$. Therefore, by \cite[4.1 and 3.6($b$)]{CEU} it follows that $\,I$ is generated by a d-sequence, hence it is of linear type (see \cite[Theorem 3.1]{dseq}). $\blacksquare$ \\
 
 An ideal $I$ for which the vanishing conditions on the Ext modules in \cref{LinTypeExt} are satisfied is called \emph{$\,(\ell-1)$-residually $S_2$}, since for $\,g \leq i \leq \ell-1\,$ every $i$-residual intersection $K$ of $I$ is such that $R/K$ satisfies Serre's condition $S_2$ (see \cite[4.1]{CEU}). We refer the interested reader to the paper of Chardin, Eisenbud and Ulrich \cite{CEU} for a thorough treatment of $s$-residually $S_2$ ideals.

 \subsection{Ideal modules}
 
 Let $R$ be a Noetherian ring. Recall that an $R$-module $E$ is called an \emph{ideal module} if $E  \neq 0$ is finitely generated and torsion-free, and moreover $E^{**}$ is free, where ${-}^*$ denotes the functor $\mathrm{Hom}_R(-, R)$.
 Equivalently, $E$ is an ideal module if and only if $E$ embeds into a finite free module $G$ with $\,\mathrm{grade}(G/E) \geq 2$ (see \cite[5.1]{SUV2003}). In particular, if $E$ is an ideal module, than $E$ has a rank, and $\mathrm{rank}_{\,}E=\mathrm{rank}_{\,}G$.
 
 An ideal module $E$ of rank 1 is isomorphic to an $R$-ideal $I$ with $\mathrm{grade}(I) \geq 2$, since one can choose $G=R$. In general, an ideal module $E$ of rank $e \geq 2$ behaves similarly to ideals of positive grade. For instance, if $\mathrm{Fitt}_e(E)$ denotes the $e$-th Fitting ideal of $E$, the non-free locus of $E$ coincides with $\, V(\mathrm{Fitt}_e(E))=\mathrm{Supp}(G/E); \,$ moreover, the analytic spread of $E$ satisfies $\, \ell(E) \geq c-e+1$, where $\,c \coloneq \mathrm{ht}(\mathrm{Fitt}_e(E))\,$ (see \cite[5.2]{SUV2003}). This is a generalization of the well-known fact that for an ideal of positive grade $\ell(I) \geq \mathrm{ht}(I)=c$. Examples of ideal modules of rank at least two include every finite direct sum of two or more ideals of grade at least two, as well as the Jacobian module of any normal complete intersection ring $\,R=k[Y_1, \ldots, Y_m]/J\,$ over a perfect field $k$.
 
 The following result shows that generic Bourbaki ideals of an ideal module have nice residual intersection properties, and will be used in the proof of \cref{thmIdealModules}.
 \begin{thm} \label{idealmodAN} \hypertarget{idealmodAN}{}
  $($\cite[5.3]{SUV2003}$)$.
 Let $R$ be a Cohen-Macaulay local ring. Every ideal module $E$ admits a generic Bourbaki ideal $I$ of height at least two. Moreover, if $E$ is free locally in codimension $s-1$, then $I$ satisfies $G_s$ and $AN_{s-1}$. 
  \end{thm}

\section{Modules of linear type} \label{SectionLinType}

 Throughout this section, we will assume the following.

\begin{set}\label{SetCMness} \hypertarget{SetCMness}{}
  \em{Let $R$ be a Noetherian local ring with $\mathrm{dim}_{\,}R=d$. Let $E$ be a finite $R$-module with $\mathrm{rank}(E)=e>0$ and analytic spread $\ell(E)=\ell$. For a fixed generating set $a_1, \ldots, a_n\,$ of $E$, let $E'$, $E''$, $\,x_1, \ldots, x_{e-1}$, $\,F'$ and $\,F''$ be constructed as in \cref{NotationBourbaki}. Assume that $E'/F'$ is isomorphic to an $R'$-ideal $J$ as an $R'$-module, let  $I =JR''$ be a generic Bourbaki ideal of $E$, and let $g= \mathrm{ht}_{\,}I$.}
\end{set}
  
 In the situation of \cref{SetCMness}, our first goal is to provide a sufficient condition for a module $E$ to be of linear type that recovers \cref{LinTypeExt} in the case of ideals. In other words, we wish to be able to deduce the linear type property of $E$ from the vanishing of finitely many modules of the form $\,\mathrm{Ext}_{R}^i(E^j, R)$, for given values of $i$ and $j$. Thanks to \cref{MainBourbaki}(c), it suffices to prove that a generic Bourbaki ideal $I$ of such a module $E$ satisfies the assumptions of \cref{LinTypeExt}.

 From the definition of generic Bourbaki ideals one has the exact sequences 
   $$\, 0 \to F' \to E' \to J \to 0 \,$$
 and 
   $$\, 0 \to F'' \to E'' \to I \to 0.\,$$
 For every $j \geq 1$, these induce $R$-epimorphisms $(E')^j \twoheadrightarrow J^j\,$ and $\, (E'')^j \twoheadrightarrow I^j$, which are obtained as the degree $j$ components of the homogeneous epimorphisms $\,\mathcal{R}(E') \twoheadrightarrow \mathcal{R}(J)\,$ and $\,\mathcal{R}(E'') \twoheadrightarrow \mathcal{R}(I)$.  Therefore, one has well-defined augmented Koszul complexes  
   $$ \, \mathbb{C}'_j \colon \, [\mathbb{K}.(x_1, \ldots, x_{e-1}; \mathcal{R}(E'))]_j \stackrel{\partial'_0}{\longrightarrow}  J^j \to 0 \,$$
 and 
   $$\, \mathbb{C}''_j \colon \, [\mathbb{K}.(x_1, \ldots, x_{e-1}; \mathcal{R}(E''))]_j \stackrel{\partial''_0}{\longrightarrow}  I^j \to 0. \,$$
 
 As we will see in several circumstances in this paper, the exactness of these complexes plays a crucial role in transferring assumptions from $E^j$ to $J^j$ or $I^j$ respectively. Our first result in this sense shows that if the complexes $\mathbb{C}'_j$ are exact locally up to a certain codimension, then the condition that $\,\mathrm{Ext}_{R}^i(E^j, R)=0\,$ for finitely many given values of $i$ and $j\,$ is preserved after replacing $E^j$ with $J^j$.

 \begin{lemma} \label{ExtPasses} \hypertarget{ExtPasses}{}
   Under the assumptions of \cref{SetCMness}, let $k$ and $s$ be positive integers with $k \leq s-2$. Assume that
   \begin{itemize}
       \item[$(i)$] $\,\mathrm{Ext}_R^{\,j+1}(E^j,R)=0$ for $1 \leq j \leq k$.       
       \item[$(ii)$] The complexes $ (\mathbb{C}'_j)_{\mathfrak{q}}$ are exact for all $\mathfrak{q} \in \mathrm{Spec}(R')$ with $\,\mathrm{depth}_{\,}R'_{\mathfrak{q}} \leq s-1$ and all $\,1 \leq j \leq k$. 
     \end{itemize}  
       Then $\,\mathrm{Ext}_{R'}^{\,j+1}(J^j,R')=0$ for $1 \leq j \leq k$.
 \end{lemma}

   \emph{Proof}. Fix $j$ with $1 \leq j \leq k$. For $i \leq j$, let $C'_i$, $Z_i$, $B_i$ and $H_i$ be the $i$th module, cycle, boundary and homology of the complex 
     $$ \,\mathbb{C}'_j \colon \, [\mathbb{K}.(x_1, \ldots, x_{e-1}; \mathcal{R}(E'))]_j \stackrel{\partial'_0}{\longrightarrow}  J^j \to 0 $$ 
   respectively. By assumption (i) we know that $ \,\mathrm{Ext}^{j+1-i}_{R'}(C'_i, R')=0\,$ for $0 \leq i \leq j$. Also, assumption (ii) implies that $\mathrm{grade}_{\,}H_i \geq s \geq k+2$. Hence, $\, \mathrm{Ext}^n_{R'}(H_i, R')=0\,$ for all $n \leq k+1$ and all $0 \leq i \leq j$. Now, by decreasing induction on $i \leq j$, we prove that $\,\mathrm{Ext}^{j+1-i}_{R'}(B_{i-1}, R')=0\,$ for $0 \leq i \leq j$. The assertion will then follow by the case $i=0$. 
   Assume that $i=j$. Since $B_j=\mathrm{im}(\partial'_{j+1})=0$ and  $\mathrm{Ext}^n_{R'}(H_j, R')=0$ for $n=0$ and $n=1$, from the exactness of $\, 0 \to B_j \to Z_j \to H_j \to 0\,$ it follows that $\, \mathrm{Hom}_{R'}(Z_j, R') \cong \mathrm{Hom}_{R'}(B_j, R')=0$. Therefore, the long exact sequence of $\mathrm{Ext}^{\cdot}_{R'}(-, R')$ induced by $\, 0 \to Z_j \to C'_j \to B_{j-1} \to 0\,$ shows that $\,\mathrm{Ext}^1_{R'}(B_{j-1}, R')=0\,$ as well. Now, assume that $j \geq i+1$ and $\,\mathrm{Ext}^{j-i}_{R'}(B_i, R')=0$. Since $\, \mathrm{Ext}^n_{R'}(H_i, R')=0\,$ for $\, n=j-i$ and $n=j-i+1$, it follows that $\mathrm{Ext}^{j-i} _{R'}(Z_i, R') \cong \mathrm{Ext}^{j-i} _{R'}(B_i, R')=0$. So, the long exact sequence of $\mathrm{Ext}^{\cdot}_{R'}(-, R')$ induced by $\, 0 \to Z_i \to C'_i \to B_{i-1} \to 0\,$ shows that $\,\mathrm{Ext}^{j+1-i}_{R'}(B_{i-1}, R')=0$, and the proof is complete. $\blacksquare$\\

 We are now ready to prove our first main result, which is a module version of \cref{LinTypeExt} and will be crucial in order to study the Cohen-Macaulay property of Rees algebras of modules in \cref{SectionCMness}.

 \begin{thm} \label{myLinType} \hypertarget{myLinType}{}
   Let $R$ be a local Gorenstein ring of dimension $d$. Let $E$ be a finite, torsion-free and orientable $R$-module with $\mathrm{rank}_{\,} E=e>0$ and $\ell(E)=\ell$. Assume that $E$ is $G_{\ell-e+2}$ and that $\,\mathrm{Ext}_R^{\,j+1}(E^j,R)= 0$ for $1 \leq j \leq \mathrm{min} \, \{\ell-e-1, \,d-3 \}$. Then, $E$ is of linear type and $E'/F'$ is isomorphic to an $R'$-ideal of linear type.
 \end{thm}

   \emph{Proof}. Without loss of generality, we may assume that $E$ is not free. Let $E=Ra_1 + \ldots + Ra_n$. Since $E$ is torsion-free, orientable and satisfies $G_{\,\ell-e+2}$, by \cref{existBourbaki}, $\,E'/F' \cong J$ and $\,E''/F'' \cong I$, where $I$ and $J$ are ideals of height at least 2, satisfying $G_{\,\ell-e+2}$, i.e. $G_{\ell(I) +1}$. Moreover, we may assume that $g \leq d-1$. In fact, if $g=d$, then $I$ is a complete intersection, hence it is of linear type by \cite[2.6]{HSV}. Then, by \cref{MainBourbaki}(c), $E$ and $J$ are of linear type. 

   If $e=1$, then $R''=R$ and $E \cong I$, an $R$-ideal of height $g$ with $2 \leq g \leq d-1$. In fact, it must be that $g=2$. Otherwise, by assumption we would have that $\mathrm{Ext}_R^{\,g-1}(I^{\,g-2},R)= 0$, contradicting the fact that $\,\mathrm{grade}_{\,}I^{\,g-2} = \mathrm{grade}_{\,}I=g$. Hence $g=2$, and so $I$ is of linear type by \cref{LinTypeExt}. 

   If $e \geq 2$, we proceed by induction on $d = \mathrm{dim}_{\,}R' \geq g \geq 2$. If $d=2$, then $g=d$ and we have already proved that $E$ and $J$ are of linear type in this case. So, assume that $d>2$. We claim that $E'_{\mathfrak{q}}$ and $J_{\mathfrak{q}}$ are of linear type for all $\,\mathfrak{q} \in \mathrm{Spec}(R')$ with $\mathrm{dim}_{\,}R'_{\mathfrak{q}} \leq d-1$. Indeed, for any such $\mathfrak{q}$ let $\mathfrak{p}= \mathfrak{q} \cap R$. Then, $E_{\mathfrak{p}}$ is a finite, torsion-free and orientable $R_{\mathfrak{p}}$-module with $\ell(E_{\mathfrak{p}}) \leq \ell$, satisfying $G_{\ell(E_{\mathfrak{p}})-e+2}$ and such that $\,\mathrm{Ext}_{R_{\mathfrak{p}}}^{\,j+1}(E_{\mathfrak{p}}^j,R_{\mathfrak{p}})= 0\,$ for $\,1 \leq j \leq \mathrm{min} \, \{\ell(E_{\mathfrak{p}})-e-1, \,\mathrm{dim}_{R_{\mathfrak{p}}}-3 \}$. Hence, by the induction hypothesis, $E_{\mathfrak{p}}$ and $J_{\mathfrak{p}} \cong E'_{\mathfrak{p}}/ F'_{\mathfrak{p}}$ are of linear type, where $E_{\mathfrak{p}}'$ and $F_{\mathfrak{p}}'$ are constructed as in \cref{NotationBourbaki} by choosing the images of $a_1, \ldots, a_n$ in $E_{\mathfrak{p}}$ as generators for $E_{\mathfrak{p}}$. Hence, their respective localizations $E'_{\mathfrak{q}}$ and $J_{\mathfrak{q}}$ are of linear type, as claimed.

   In particular, by \cref{SUV3.11}, the complexes $\, (\mathbb{C}'_j)_{\mathfrak{q}}$ are exact for all $j$ and all $\,\mathfrak{q} \in \mathrm{Spec}(R')\,$ with $\,\mathrm{depth}_{\,}R'_{\mathfrak{q}} \leq d-1$. Therefore, by \cref{ExtPasses}, $\,\mathrm{Ext}_{R'}^{\,j+1}(J^j,R')=0\,$ for $\, 1 \leq j \leq\mathrm{min} \, \{\ell-e-1, \,d-3 \}$. Hence, for $j$ in the same range, also $\,\mathrm{Ext}_{R''}^{\,j+1}(I^j,R'')=0$. Now, if $\,3 \leq g \leq d-1$, we could choose $\,j=g-2\,$ to obtain $\,\mathrm{Ext}_{R''}^{\,g-1}(I^{g-2},R'')=0$. But this contradicts the fact that $\mathrm{grade}_{\,}I^{g-2}=\mathrm{grade}_{\,}I=g$. So, it must be that $g=2$, whence $I$ is of linear type by \cref{LinTypeExt}. $\blacksquare$\\

 In the remaining part of this section we show that the assumptions of \cref{myLinType} are satisfied by several classes of modules. Our first example is given by modules of projective dimension one, improving the well-known result that a module of projective dimension one is of linear type when it satisfies $G_{\infty}$ (see \cite[Propositions 3 and 4]{Avramov}, \cite[1.1]{HuSym} and \cite[3.4]{SVsym}).

 \begin{prop} \label{projdim1} \hypertarget{projdim1}{}
   Let $R$ be a local Cohen-Macaulay ring, $E$ a finite, torsion-free $R$-module with $\mathrm{rank}_{\,} E=e>0$ and $\ell(E)=\ell$. Assume that $\mathrm{projdim}_{\,}E=1$ and that $E$ satisfies $G_{\ell-e+2}$. Then, $E$ has no proper reductions and $\,\mathrm{Ext}_R^{\,j+1}(E^j,R)= 0$ for all $j \geq 1$. 
 \end{prop}

   \emph{Proof}. Without loss of generality, we may assume that the residue field of $R$ is infinite. Let $\mathfrak{p} \in \mathrm{Spec}(R)$ be such that $\,\mathrm{dim}_{\,}R_{\mathfrak{p}} \leq \ell -e +1$. Then, $E_{\mathfrak{p}}$ satisfies $\,G_{\infty}$ and has projective dimension one. Hence, by \cite[Proposition 3 and 4]{Avramov} it follows that $\,\mathcal{R}(E_{\mathfrak{p}}) \cong \mathcal{S}(E_{\mathfrak{p}})\,$ is a complete intersection. In particular, $\,E_{\mathfrak{p}}$ is of linear type and $\,\mathrm{projdim}_{\,}E_{\mathfrak{p}}^j \geq j\,$ for all $\,j \geq 1$. Hence, $\,\mathrm{Ext}_{R_{\mathfrak{p}}}^{j+1} (E_{\mathfrak{p}}, R_{\mathfrak{p}})=0\,$ for all $j \geq 1$. 

   Now, let $U$ be a minimal reduction of $E$, generated by $\ell$ general elements. Since for all $\,\mathfrak{p} \in \mathrm{Spec}(R)\,$ with $\,\mathrm{dim}_{\,}R_{\mathfrak{p}} \leq \ell -e +1\,$ $\,E_{\mathfrak{p}}$ is of linear type, then for any such $\mathfrak{p}$ one has $(E/U)_{\mathfrak{p}}=0$. Hence, $\,\mathrm{ht}(F_{0}(E/U)) \geq \ell - e +2, \,$ where $\,F_0(E/U)\,$ is the $0$th Fitting ideal of $E/U$. On the other hand, since $E/U$ is generated by $\,\mu(E)-\ell\,$ elements, by Eagon-Northcott's Theorem \cite[Theorem 3]{EN}, if $\,F_0(E/U) \neq R\,$ we would have that $\,\mathrm{ht}(F_0(E/U)) \leq \ell - e +1$. This is a contradiction, hence it must be that $E=U$, in which case $E$ has no proper reductions and $\,\mu(E)= \ell$. 
   
   In fact, since $\,\ell(E_{\mathfrak{p}}) \leq \ell \,$ for every $\mathfrak{p} \in \mathrm{Spec}(R)$, one can repeat the same argument for a minimal reduction $U_{\mathfrak{p}}$ of $E_{\mathfrak{p}}$ to show that $\,\mu(E_{\mathfrak{p}})= \ell(E_{\mathfrak{p}}) \leq \mathrm{dim} R_{\mathfrak{p}} -e+1\,$ for all $\mathfrak{p} \in \mathrm{Spec}(R)$. In particular, $E$ satisfies $\,G_{\infty}$, and the argument at the beginning of the proof then finally proves that $\,\mathrm{Ext}_R^{\,j+1}(E^j,R)= 0$ for all $j \geq 1$. $\blacksquare$\\
   
 Another class of modules satisfying \cref{myLinType} can be constructed using strongly Cohen-Macaulay ideals. Recall that an ideal $I$ in a local Cohen-Macaulay ring $R$ is \emph{strongly Cohen-Macaulay} if all of its Koszul homologies are Cohen-Macaulay $R$-modules \cite{SCM}. 

 \begin{prop} \label{SCM} \hypertarget{SCM}{}
   Let $R$ be a local Gorenstein ring, with $\,\mathrm{dim}_{\,}R=d$. Let $I$ be a strongly Cohen-Macaulay ideal of height two with $\,\ell(I)=\ell$, satisfying $\,G_{\ell +1}$. Let $F$ be a free $R$-module of rank $\,e-1 >0$, and let $\,E= I \oplus F$. Then $\,\ell(E) = \ell +e -1, \,$ $E$ satisfies $\,G_{\ell(E) +e-2}$ and $\,\mathrm{Ext}^{j+1}_R (E^j, R)=0\,$ for $\,1 \leq j \leq \ell(E)-e+1$.
 \end{prop}

   \emph{Proof}. Notice that $E= I \oplus F$ is a torsion free $R$-module with $\mathrm{rank}E = e$ and satisfies $G_{\ell +1}$ since $I$ does. Moreover, since $\mathcal{S}(F)$ is a polynomial ring in $e-1$ variables over $R$, the special fiber ring $\,\mathcal{F}(E) \cong \mathcal{F}(I) \otimes_R \mathcal{S}(F)\,$ has dimension $\ell(E)= \ell +e -1$. In particular, $E$ satisfies $\,G_{\ell(E) -e +2}$.

   Now, since $I$ is strongly Cohen-Macaulay of height 2 and satisfies $G_{\ell +1}$, by \cite[2.10]{Ulrich} we have that $\,\mathrm{depth}_{\,}I^j \geq d-j\,$ for $\,1 \leq j \leq \ell$. Since $R$ is Gorenstein, then $\,\mathrm{Ext}^{j+1}_R (I^j, R)=0 \,$ for $\,1 \leq j \leq \ell$, and thus $I$ is of linear type by \cref{LinTypeExt}. In particular, for all $j$ 
      $$ \,[\mathcal{S}(E)]_j \cong [\mathcal{R}(I) \otimes_R \mathcal{S}(F)]_j = \bigoplus_{0 \leq i \leq j\,} (I^i \otimes_R F^{j-i}) \, $$
   is a torsion free $R$-module. Therefore, $ \, \displaystyle{E^j \cong [\mathcal{S}(E)]_j \cong \bigoplus_{0 \leq i \leq j \,} (I^i \otimes_R F^{j-i})},\,$ whence $\,\mathrm{Ext}^{j+1}_R (E^j, R)=0\,$ for $1 \leq j \leq \ell= \ell(E)-e+1$.  $\blacksquare$\\

 We now provide a general technique to construct modules of linear type, as submodules of modules that satisfy \cref{myLinType}. 
 
 \begin{thm} \label{LinTypeSub}
   Let $R$ be a local Gorenstein ring. Let $M$ be a finite, torsion-free $R$-module with $\,\mathrm{rank}_{\,} M=e>0\,$ and $\,\ell(M)=\ell$. Assume that $\,\mathrm{Ext}_R^{\,j+1}(M^j,R)= 0\,$ for $\,1 \leq j \leq \mathrm{min} \, \{\ell-e-1, \,d-3 \}\,$ and that $M$ satisfies $\,G_{\ell-e+2}$. 
    
   Let $E$ be an orientable submodule of $M$ with $\, \mathrm{dim}_{\,}(M/E) \leq \mathrm{max}_{\,} \{\, d- \ell +e -2, \,0 \, \}$. If $\,\ell=d+e-1,\,$ assume also that $E$ is generated by $\,\ell$ elements. Then $E$ is of linear type.
 \end{thm}

   \emph{Proof}. We may assume that $d>0$. We show that $E$ satisfies the assumptions of \cref{myLinType}. First, notice that since $\mathrm{dim}(M/E) \leq \mathrm{max}_{\,} \{\, d- \ell +e -2, \,0 \, \},\,$ it follows that $\,E_{\mathfrak{p}} = M_{\mathfrak{p}}$ for every prime ideal $\mathfrak{p}$ with $\,\mathrm{dim}_{\,}R_{\mathfrak{p}} < \mathrm{min}_{\,} \{\, \ell -e +2, d\, \}$. In particular, $E$ has a rank, and $\mathrm{rank}_{\,}E=e$. 
 
   Moreover, $E$ is $G_{\ell - e+2}$. This is clear from the discussion above if $\,\ell \leq d+e-2$. Otherwise, by \cref{anspread} one must have that $\,\ell = d+e-1$, in which case the discussion above shows that $E$ satisfies $G_d$. Since by assumption $\,\mu(E) \leq \ell$, it follows that actually $E$ is $G_{\,\ell-e+2}$, as claimed. Moreover, since $\,\ell(E) \leq \mu(E) \leq \ell$, then $E$ satisfies $G_{\ell(E) - e+2}$.
 
   Therefore, it only remains to show that our assumptions imply the vanishing of the Ext modules in the statement of \cref{myLinType}. Since $\,E_{\mathfrak{p}} = M_{\mathfrak{p}}$ for every prime ideal $\mathfrak{p}$ with $\,\mathrm{dim}_{\,}R_{\mathfrak{p}} < \mathrm{min}_{\,} \{\, \ell -e +2, d\, \}$, it follows that for all $j \geq 1$ $\mathrm{grade}(M^j /E^j) \geq \mathrm{min}_{\,} \{\, \ell -e +2, d\, \}$. Since $R$ is Gorenstein, this implies that $\,\mathrm{Ext}^{j+2}_R(M^j /E^j,R)= 0\,$ for $\,1 \leq j \leq \mathrm{min} \, \{\ell-e-1, \,d-3 \}$. The long exact sequence of $\,\mathrm{Ext}^{\bullet}_R(-,R)\,$ induced by the exact sequence
    $$ \, 0 \to E^j \to M^j \to M^j/E^j \to 0 \, $$
   then shows that $\,\mathrm{Ext}_R^{\,j+1}(E^j,R)= 0$ for $1 \leq j \leq \mathrm{min} \, \{\ell-e-1, \,d-3 \}$. Hence, $E$ is of linear type by \cref{myLinType}. $\blacksquare$\\
 
 Explicit examples of submodules $E$ of $M$ satisfying the assumptions of \cref{LinTypeSub} can be constructed by means of the following lemma.
 
 \begin{lemma} \label{generalBourbaki} \hypertarget{generalBourbaki}{}
   Let $(R, \mathfrak{m})$ be a local Cohen-Macaulay ring with infinite residue field, let $M$ be a finitely generated, torsion-free $R$-module with $\mathrm{rank}_{\,}M=e>0$, satisfying $G_s$ for some $s\geq 2$. Let $E$ be a submodule generated by $\,s+e-1\,$ general elements in $\mathfrak{m}M$. Then, $\,\mathrm{dim}(M/E) \leq d-s$.
 \end{lemma}
 
 \emph{Proof}. Since $E$ is generated by $\,s+e-1\,$ general elements $\,x_1, \ldots, x_{s+e-1}\,$ in $\mathfrak{m}M$, it makes sense to define the module $\,\overline{M} \coloneq M/(Rx_1+ \ldots+ Rx_{e-1})$. Since $M$ is torsion-free and satisfies $G_s$ with $s\geq 2$, by \cite[2.2(a) and (b)]{CPU2003} it follows that $\overline{M}$ is a torsion-free module of rank 1 and satisfies $G_s$. That is, $\overline{M} \cong I$, an $R$-ideal satisfying $G_s$. 
 
 Now, let $J$ be the ideal generated by the images of $\,x_e, \ldots, x_{s+e-1}\,$ in $I$. Since the $\,x_e, \ldots, x_{s+e-1}\,$ are general elements, by \cite[1.6]{Ulrich} it follows that $\,J \colon_{\!R\,}I\,$ is an $s$-residual intersection of $I$, whence $\, \mathrm{ht}(J \colon_{\!R\,}I) \geq s$. Hence, 
  $$ \,\mathrm{dim}(M/E) = \mathrm{dim}_{\,}(I/J) =\mathrm{dim}_{\,}R-\mathrm{ht}(J \colon_{\!R\,}I) \leq \mathrm{dim}_{\,}R-s, \,$$
 as we wanted to prove. $\blacksquare$\\

 \begin{cor}
   Let $R$ be a local Gorenstein UFD, of dimension $d$ and infinite residue field. Let $M$ be a finite, torsion-free $R$-module with $\mathrm{rank}_{\,} M=e>0$ and $\ell(M)=\ell$. Assume that $\,\mathrm{Ext}_R^{\,j+1}(M^j,R)= 0$ for $1 \leq j \leq \mathrm{min} \, \{\ell-e-1, \,d-3 \}\,$ and that $M$ satisfies $G_{\ell-e+2}$. 
   
   If $\,\ell \leq d+e-2$, let $E$ be a submodule generated by $\,\ell+1\,$ general elements inside $\, \mathfrak{m}M$, while if $\,\ell=d+e-1, \,$ let $E$ be a submodule generated by $\,\ell$ general elements inside $\, \mathfrak{m}M$. Then, $E$ is of linear type.
 \end{cor}
 
   \emph{Proof}. Since $M$ satisfies $G_{\ell-e+2}$, from  \cref{generalBourbaki} it follows that
    $$ \,\mathrm{dim}(M/E) \leq \Big\lbrace \begin{array}{cc}
                       d-\ell+e-2\, \; \; & \mathrm{if} \; \; \ell \leq d+e-2 \\
                       0 \; & \mathrm{if} \; \; \ell=d+e-1 \\
                   \end{array} $$
   Moreover, $E$ is orientable since $R$ is assumed to be a UFD. Therefore, the conclusion follows from \cref{LinTypeSub}. $\blacksquare$\\
 
 Combining the latter result with \cref{projdim1} and \cref{SCM} respectively, we have the following corollaries.
 
 \begin{cor}
   Let $R$ be a local Gorenstein UFD, of dimension $d$ and infinite residue field. Let $M$ be a finite $R$-module with $\mathrm{rank}_{\,} M=e>0$,  $\,\ell(M)=\ell\,$ and $\, \mathrm{projdim}_{\,}M=1,\,$ satisfying $G_{\ell-e+2}$. 
   
   If $\,\ell \leq d+e-2, \,$ let $E$ be a submodule generated by $\,\ell+1\,$ general elements inside $\, \mathfrak{m}M$, while if $\,\ell=d+e-1, \,$ let $E$ be a submodule generated by $\,\ell$ general elements inside $\, \mathfrak{m}M$. Then, $E$ is of linear type.
 \end{cor}

 \begin{cor}
   Let $R$ be a local Gorenstein UFD, of dimension $d$ and infinite residue field. Let $I$ be a strongly Cohen-Macaulay ideal of height two with $\ell(I)=\ell$ satisfying $G_{\ell +1}$, $F$ a free $R$-module of rank $e-1 >0$, and let $M= I \oplus F$. 
   
   If $\,\ell \leq d+e-2, \,$ let $E$ be a submodule generated by $\,\ell+1\,$ general elements inside $\, \mathfrak{m}M$, while if $\,\ell=d+e-1,\,$ let $E$ be a submodule generated by $\,\ell$ general elements inside $\, \mathfrak{m}M$. Then, $E$ is of linear type.
 \end{cor}

  Finally, notice that since $R$ is Gorenstein the assumption that $\,\mathrm{Ext}_R^{\,j+1}(E^j,R)= 0\,$ for $\,1 \leq j \leq \mathrm{min} \, \{\ell-e-1, \,d-3 \}\,$ in \cref{myLinType} is always satisfied if $\,\mathrm{depth}_{\,}E^j \geq d-j\,$ for $\,1 \leq j \leq \ell-e-1,\,$ a condition that reminds of \cref{LinTypeDepth}. In fact, \cref{herLinType} at the end of the next section will provide a  stronger result, which extends \cref{LinTypeDepth} to modules. In order to prove \cref{herLinType}, we first need to understand how lower bounds on the depths of powers of a module behave under passage to a generic Bourbaki ideal, which we will discuss in \cref{SectionCMness}.

\section{Modules with Cohen-Macaulay Rees algebra} \label{SectionCMness}

 In this section we study the Cohen-Macaulay property of Rees algebras of modules. Our main goal is to show that the Rees algebra $\mathcal{R}(E)$ of a module $E$ is Cohen-Macaulay if finitely many powers of $E$ have sufficiently high depth, so that we recover \cref{JUGNN} in the case when $E$ is an ideal. 
 
 With the assumptions and notations of \cref{SetCMness}, the first step is to show that lower bounds on the depths of finitely many powers $E^j$ are preserved when passing to $I^j$, where $I$ is a generic Bourbaki ideal of $E$. Similarly as in the proof of \cref{myLinType}, the exactness of the complexes
   $$ \, \mathbb{C}'_j \colon \, [\mathbb{K}.(x_1, \ldots, x_{e-1}; \mathcal{R}(E'))]_j \stackrel{\partial'_0}{\longrightarrow}  J^j \to 0 \, $$
 and 
    $$ \, \mathbb{C}''_j \colon \, [\mathbb{K}.(x_1, \ldots, x_{e-1}; \mathcal{R}(E''))]_j \stackrel{\partial''_0}{\longrightarrow}  I^j \to 0 \, $$
 described in \cref{SectionLinType} will be crucial, as explained in the following lemma.
      
 \begin{lemma} \label{DepthPasses} 
    Under the assumptions of \cref{SetCMness}, let $k$ and $s$ be integers such that $\,\mathrm{depth}_{\,}E^j \geq s-j\,$ for $\,0 \leq j \leq k$. If the complexes $\mathbb{C}''_j$ are exact for $\,0 \leq j \leq k$, then $\,\mathrm{depth}_{\,} I^j \geq s-j\,$ for $\,0 \leq j \leq k$.
 \end{lemma}

  \emph{Proof}. Fix $j$ with $\,0 \leq j \leq k$. For $\,0 \leq i \leq j, \,$ let $C''_i$ and $B_i$ be the $i$th module and boundary of the complex $\, \displaystyle \mathbb{C}''_j$ respectively. From our assumptions it follows that $\,\mathrm{depth}_{\,} (E'')^j \geq s-j\,$ for $\,0 \leq j \leq k.\,$ Hence, by the Depth Lemma, for $\, 0 \leq i \leq j\,$ one has that  $\,\mathrm{depth}_{\,}C''_i \geq s-j+i$. Now this implies that $\,\mathrm{depth}_{\,}B_{i-1} \geq s-j+i\,$ for all $0 \leq \,i \leq j$, which can be proved by decreasing induction on $i$, using again the Depth Lemma. $\blacksquare$\\

 The next lemma clarifies how the complexes $\mathbb{C}'_j$ may fail to be exact.

 \begin{lemma} \label{SmallCodim}
     Under the assumptions of \cref{SetCMness}, let $k$ and $s$ be integers such that $0 \leq k \leq s$. Assume that $\,\mathrm{depth}_{\,}E^j \geq d-s+k-j\,$ for $\,0 \leq j \leq k$. Then, for all $\,0 \leq j \leq k$ and every minimal prime $\mathfrak{q}$ in $\,\mathrm{Supp} \Big( \bigoplus_{i=0}^j H_i(\mathbb{C}'_j) \Big)$, $\,\mathrm{dim}_{\,}R'_{\mathfrak{q}} \leq s$.
 \end{lemma}

  \emph{Proof}. Fix $j$ with $\,0 \leq j \leq k$ and suppose that $\,\mathrm{dim}_{\,}R'_{\mathfrak{q}} \geq s+1\,$ for some prime ideal $\mathfrak{q}$ as in the statement. We will prove that $(\mathbb{C}'_j)_{\mathfrak{q}}$ is exact, and this will contradict the fact that $H_i(\mathbb{C}'_j)_{\mathfrak{q}} \neq 0\,$ for some $i$ with $\,0 \leq i \leq j$. 

  First, notice that $H_i(\mathbb{C}'_j)_{\mathfrak{q}}$ is zero or has depth zero for all $0 \leq i \leq j$, since $\mathfrak{q}$ is minimal in $\mathrm{Supp} \Big( \bigoplus_{i=0}^j H_i(\mathbb{C}'_j) \Big)$. Now, let $\mathfrak{p}= \mathfrak{q} \cap R$. Since $R_{\mathfrak{p}} \to R'_{\mathfrak{q}}$ is a flat local map, then $\,\mathrm{dim}_{\,}R'_{\mathfrak{q}} -\mathrm{depth}(E'_{\mathfrak{q}})^j = \mathrm{dim}_{\,}R_{\mathfrak{p}} -\mathrm{depth}_{\,}(E_{\mathfrak{p}})^{\,j} \leq d - \mathrm{depth}_{\,}E^j\,$ \cite[Theorem 15.1, Theorem 23.3 and Exercise 17.5]{Mats}. Hence, from our assumption on $\mathrm{depth}_{\,}E^{j-i} $ it follows that $\,\mathrm{depth}(E'_{\mathfrak{q}})^{j-i} \geq \mathrm{dim}_{\,}R'_{\mathfrak{q}}-s+k-j+i \geq i+1$. Hence, $(\mathbb{C}'_j)_{\mathfrak{q}}$ is acyclic by the Acyclicity Lemma, and therefore it is exact, since $H_{-1}(\mathbb{C}'_j)=0$ by construction. $\blacksquare$\\

 We are now ready to state and prove the desired generalization of \cref{JUGNN}. 
    
 \begin{thm} \label{GenGNN} \hypertarget{GenGNN}{}
    Let $R$ be a local Gorenstein ring of dimension $d$ with infinite residue field. Let $E$ be a finite, torsion-free, orientable $R$-module, with $\mathrm{rank}E=e>0$ and $\ell(E)=\ell$. Let $g$ be the height of a generic Bourbaki ideal of $E$, and assume that the following conditions hold.
    \begin{itemize}
       \item[$(a)$] $E$ satisfies $G_{\,\ell-e+1}$. 
       \item[$(b)$] $r(E) \leq k$ for some integer $1 \leq k \leq \ell-e$. 
       \item[$(c)$] $\displaystyle{
             \mathrm{depth}_{\,}E^j \geq \Big\lbrace \begin{array}{cc}
                       d-g-j+2\, \; \; & \mathrm{for} \; \; 1 \leq j \leq \ell-e-k-g+1 \\
                       d-\ell+e+k-j\, \; & \mathrm{for} \; \; \ell-e-k-g+2 \leq j \leq k \\
                   \end{array}}$
       \item[$(d)$] If $g = 2$, $\,\mathrm{Ext}_{R_{\mathfrak{p}}}^{\,j+1}(E_{\mathfrak{p}}^j , R_{\mathfrak{p}}) =0\,$ for $\,\ell-e-k \leq j \leq \ell-e-3\,$ and for all $\,\mathfrak{p} \in \mathrm{Spec}(R)$ with $\, \mathrm{dim} R_{\mathfrak{p}} = \ell-e\, $ such that $E_p$ is not free.
      \end{itemize}
    Then, $\mathcal{R}(E)$ is Cohen-Macaulay.
 \end{thm}

   \emph{Proof}. Let $E=Ra_1 + \ldots + Ra_n$. By \cref{existBourbaki}, $E$ admits a generic Bourbaki ideal $I$ with $\mathrm{ht}I = g \geq 2$ and $r(I) \leq r(E) \leq k$, that satisfies $G_{\,\ell-e+1}$, i.e. $G_{\,\ell(I)}$. We will prove that $\mathcal{R}(I)$ is Cohen-Macaulay, so that $\mathcal{R}(E)$ is Cohen-Macaulay by \cref{MainBourbaki}. 

   If $e=1$, then $R''=R$ and $E \cong I$, an $R$-ideal with $\ell(I)=\ell$ and such that
     $$ \, \mathrm{depth}_{\,}E^j \geq \Big\lbrace \begin{array}{cc}
                       d-g-j+2\, \; \; & \mathrm{for} \; \; 1 \leq j \leq \ell-k-g \\
                       d-\ell+k-j+1\, \; & \mathrm{for} \; \; \ell-k-g+1 \leq j \leq k \\
                   \end{array} \, $$
   Since $R$ is Gorenstein, by \cref{CreateAN} it follows that $I$ satisfies $AN_{\ell-k-1}$. Moreover, all assumptions of \cref{JUGNN} are satisfied, hence $\mathcal{R}(I)$ is Cohen-Macaulay. 

   Now, assume $e \geq 2$. It suffices to prove that 
     $$\,\mathrm{depth}_{\,}I^j \geq \Big\lbrace \begin{array}{cc}
                       d-g-j+2\, \; \; & \mathrm{for} \; \; 1 \leq j \leq \ell-e-k-g+1 \\
                       d-\ell+e+k-j\, \; & \mathrm{for} \; \; \ell-e-k-g+2 \leq j \leq k \\
                   \end{array}\, $$
                   
   In fact, then $\,\mathrm{depth}_{\,}I^j \geq d-\ell(I)+k-j+1\,$ for $\,1 \leq j \leq k\,$ and by \cref{CreateAN} $I$ would satisfy $AN_{\ell-k-e}$, i.e. $AN_{\ell(I)-k-1}$. Hence, $\mathcal{R}(I)$ would be Cohen-Macaulay by \cref{JUGNN}. 
 
   Notice that $\,\mathrm{depth}_{\,}(E'')^j \geq \mathrm{depth}_{\,} E^j$, since $R''$ is flat over $R$. Also, if $g \geq 3$, then by \cref{existBourbaki}(c) $I$ is a free direct summand of $E''$, so for every $j$, $I^j$ satisfies the same depth condition as $(E'')^j$, hence as $E^j$. Thus, we may assume that $g=2$, so that assumption (c) becomes 
     $$ \,\mathrm{depth}_{\,}E^j \geq \Big\lbrace \begin{array}{cc}
                       d-j\, \; \; & \quad \; \; \mathrm{for} \; \; 1 \leq j \leq \ell-e-k-1 \\
                       d-\ell+e+k-j\, \; & \mathrm{for} \; \; \ell-e-k \leq j \leq k \\
                   \end{array}\, $$          
   In particular, $\mathrm{depth}_{\,}E^j \geq d-\ell+e+k-j$ for $1\leq j \leq k$, and also for $j=0$ since $R$ is Cohen-Macaulay. Therefore, by \cref{DepthPasses}  (with $s=d-\ell+e+k$) and \cref{SmallCodim} (with $s=\ell-e$), it suffices to show that the complexes $(\mathbb{C}'_j)_{\mathfrak{q}}$ are exact for all $\mathfrak{q} \in \mathrm{Spec}(R')$ with $\,\mathrm{dim}_{\,}R'_{\mathfrak{q}} \leq \ell-e\,$ and all $0 \leq j \leq k$. 

   For each such $\mathfrak{q}$, let $\mathfrak{p}=\mathfrak{q} \cap R$. If $E_{\mathfrak{p}}$ is free, then also its localization $E'_{\mathfrak{q}}$ is free. Hence, $E'_{\mathfrak{q}}$ is of linear type and $\mathcal{R}(E'_{\mathfrak{q}})$ is Cohen-Macaulay, so by \cref{MainBourbaki}(c) and \cref{SUV3.11}, the complexes $(\mathbb{C}'_j)_{\mathfrak{q}}$ are exact for all $j$. If $E_{\mathfrak{p}}$ is not free, then by assumption $\,\mathrm{Ext}_{R_{\mathfrak{p}}}^{\,j+1}(E_{\mathfrak{p}}^{\,j} , R_{\mathfrak{p}}) =0\,$ for all $\,\mathfrak{p} \in \mathrm{Spec}(R)$ with $\, \mathrm{dim} R_{\mathfrak{p}} \leq \ell-e\,$ and all $j$ with $\,\ell-e-k \leq j \leq \ell-e-3$. In particular, since $\,\ell(E_{\mathfrak{p}}) \leq \ell \,$ and $\,\mathrm{dim}_{\,}R'_{\mathfrak{q}} \leq \ell-e,\,$ this is true for $\,\ell-e-k \leq j \leq \mathrm{min} \, \{\ell(E_{\mathfrak{p}})-e-1, \,\mathrm{dim}_{\,}R_{\mathfrak{p}}-3 \}$. Moreover, since $R_{\mathfrak{p}}$ is Gorenstein, by assumption (c) the same vanishing holds if $\,1 \leq j \leq \ell-e-k-1$. Hence, $\,\mathrm{Ext}_{R_{\mathfrak{p}}}^{\,j+1}(E_{\mathfrak{p}}^{\,j} , R_{\mathfrak{p}}) =0\,$ for all $\,1 \leq j \leq \mathrm{min} \, \{\ell(E_{\mathfrak{p}})-e-1, \,\mathrm{dim}_{\,}R_{\mathfrak{p}}-3 \}$. Therefore, by \cref{myLinType}, $\,E_{\mathfrak{p}}'/ F_{\mathfrak{p}}'\,$ is isomorphic to an $R_{\mathfrak{p}}'$-ideal of linear type, where $E_{\mathfrak{p}}'$ and $F_{\mathfrak{p}}'$ are constructed as in \cref{NotationBourbaki} by choosing the images of $a_1, \ldots, a_n$ in $E_{\mathfrak{p}}$ as generators for $E_{\mathfrak{p}}$. Therefore, $\,(E'/ F')_{\mathfrak{q}}\,$, which is a localization of $\,E_{\mathfrak{p}}'/ F_{\mathfrak{p}}',\,$ is isomorphic to an $R'_{\mathfrak{q}}$-ideal of linear type. Thus, all the $(\mathbb{C}'_j)_{\mathfrak{q}}$ are exact by \cref{SUV3.11}. $\blacksquare$\\

 Comparing the statements of \cref{JUGNN} and of \cref{GenGNN}, it is clear that assumption (d) in \cref{GenGNN} is redundant in case $E$ is of rank one. Moreover, in this case assumption (c) simplifies to 
   $$\, (c') \quad \mathrm{depth}_{\,}E^j \geq d-\ell+e+k-j \quad \mathrm{for} \quad 1\leq j \leq k.$$
 It is then natural to ask what other situations produce similar simplified statements for \cref{GenGNN}. 

 Notice that assumption (d) in \cref{GenGNN} is vacuously satisfied when $\,k \leq 2,\,$ or when $\,\ell -e \leq 3$. Also, if $\,\ell -e \leq 2$, assumption (c) can be replaced by assumption (c') above. In particular, we have the following two corollaries, which correspond to the cases when a generic Bourbaki ideal of $E$ has analytic deviation at most $1$ or at most $2$, respectively.

  \begin{cor} 
    Let $R$ be a local Gorenstein ring of dimension $d$ with infinite residue field. Let $E$ be a finite, torsion-free, orientable $R$-module, with $\mathrm{rank}E=e>0$ and $\ell(E)=\ell$. Let $g$ be the height of a generic Bourbaki ideal of $E$, and assume that the following conditions hold.
    \begin{itemize}
       \item[$($a$)$] $E$ satisfies $G_{\,\ell-e+1}$. 
       \item[$($b$)$] $r(E) \leq k$ for some integer $\,1 \leq k \leq \ell-e\,$ and $\,\ell -e \leq 2$.
       \item[$($c$)$] $\mathrm{depth}_{\,}E^j \geq d-\ell+e+k-j\,$ for $\,1 \leq j \leq k$.
      \end{itemize}
    Then, $\mathcal{R}(E)$ is Cohen-Macaulay.
  \end{cor} 

  \begin{cor} 
    Let $R$ be a local Gorenstein ring of dimension $d$ with infinite residue field. Let $E$ be a finite, torsion-free, orientable $R$-module, with $\mathrm{rank}E=e>0$ and $\ell(E)=\ell$. Let $g$ be the height of a generic Bourbaki ideal of $E$, and assume that the following conditions hold.
    \begin{itemize}
       \item[$($a$)$] $E$ satisfies $G_{\,\ell-e+1}$. 
       \item[$($b$)$] $r(E) \leq k$ for some integer $\,1 \leq k \leq \ell-e\,$ and $\,\ell -e \leq 3$.
       \item[$($c$)$] $\mathrm{depth}_{\,}E \geq d-g+1$ and $\,\mathrm{depth}_{\,}E^j \geq d-\ell+e+k-j\,$ for $\,2 \leq j \leq k$.
      \end{itemize}
    Then, $\mathcal{R}(E)$ is Cohen-Macaulay.
  \end{cor}

  Observe also that assumption (b) of \cref{GenGNN} implies that $\,r(E) \leq \ell(E)-e, \,$ which is a necessary condition in order for $\,\mathcal{R}(E)\,$ to be Cohen-Macaulay \cite[4.2]{SUV2003}. Recalling that $\,\ell(E)-e \leq d-1,\,$ it then makes sense to investigate the Cohen-Macaulay property of $\,\mathcal{R}(E)\,$ for rings of small dimensions. In particular, in \cite[4.4 and 4.6]{SUV2003}, sufficient conditions were given in the case when $d \leq 5$, under the assumption that $r(E) \leq 2$ and $E$ satisfies $G_d$. As a consequence of \cref{GenGNN}, when $d \leq 5$ we can now give similar results for the full range of non-zero admissible reduction numbers, assuming that $E$ satisfies the weaker condition $G_{\ell-e+1}$. 

  \begin{cor} 
    Let $R$ be a Gorenstein local ring of dimension $4$ with infinite residue field. Let $E$ be a finite, torsion-free, orientable $R$-module with $\mathrm{rank}_{\,}E =e$ and $\ell(E)=\ell$. Assume that $E$ satisfies $G_{\ell-e+1}$ and $\,r(E) \leq \ell-e$. Then, $\mathcal{R}(E)$ is Cohen-Macaulay if one of the following conditions holds.
    \begin{itemize}
       \item[$($a$)$] $r(E)=1 < \ell-e\,$ and $\, \mathrm{depth}_{\,}E \geq 2$, or $\, r(E)=1=\ell-e\,$ and $\,\mathrm{depth}_{\,}E \geq 3$.
       \item [$($b$)$] $r(E)=2 < \ell-e\,$ and $\,\mathrm{depth}_{\,}E \geq 2$, or $\, r(E)=2=\ell-e\,$ and $\,\mathrm{depth}_{\,}E^j \geq 4-j\,$ for $\,1 \leq j \leq 2$.
       \item[$($c$)$] $r(E)=3\,$  and $\,\mathrm{depth}_{\,}E^j \geq 4-j\,$ for $\,1 \leq j \leq 3$.
      \end{itemize}
  \end{cor}

  \emph{Proof}. Since $\,\ell-e  \leq d-1=3, \,$ if  $\,\ell-e=3\,$ it follows that $E$ satisfies $G_d$. Hence, if $r(E) \leq 2$, by \cite[4.6(b)]{SUV2003} it suffices to assume that $\, \mathrm{depth}_{\,}E \geq 2$. In all the remaining cases, our assumptions imply the assumptions of \cref{GenGNN} with $d=4$ and $k=r(E)$. In fact, since $\,r(E) \leq \ell -e \leq 3, \,$ assumption (d) in \cref{GenGNN} is vacuously satisfied, and if $\,r(E)=3\,$ then it must be $\,\ell - e=3$. $\blacksquare$

 \begin{cor} 
    Let $R$ be a Gorenstein local ring of dimension $5$ with infinite residue field. Let $E$ be a finite, torsion-free, orientable $R$-module with $\mathrm{rank}_{\,}E =e$ and $\ell(E)=\ell$. Assume that $E$ satisfies $G_{\ell-e+1}$. Let $g$ be the height of a generic Bourbaki ideal of $E$. Then, $\mathcal{R}(E)$ is Cohen-Macaulay if one of the following conditions holds.    
     \begin{itemize}
       \item[$($a$)$] $\,\ell-e =4, \,$ $\,r(E) \leq 2\,$ and $\,\mathrm{depth}_{\,}E \geq 4$.
       \item[$($b$)$] $\,\ell-e =4$, $\,r(E) \geq \ell-e -1=3$, and $\,\mathrm{depth}_{\,}E^j \geq r(E) +1 -j\,$ for $\,1 \leq j \leq r(E)$. If $g=2, \,$ assume also that $\,\mathrm{Ext}^2_{R_{\mathfrak{p}}} (E_{\mathfrak{p}}, R_{\mathfrak{p}})=0\,$ for all $\,\mathfrak{p} \in \mathrm{Spec}(R)$ with $\, \mathrm{dim} R_{\mathfrak{p}} = 4\,$ such that $E_p$ is not free. 
       \item[$($c$)$] $\,r(E)= \ell-e \leq 3, \,$ and $\,\mathrm{depth}_{\,}E^j \geq 5-j\,$ for $\,1 \leq j \leq r(E)$.
       \item[$($d$)$] $\,r(E)= \ell-e -1 \leq 2, \,$ and $\,\mathrm{depth}_{\,}E^j \geq 4-j\,$ for $\,1 \leq j \leq r(E)$.     
       \item[$($e$)$] $\,\ell-e =3, \,$ $\,r(E)=1$, and $\displaystyle{\,
             \mathrm{depth}_{\,}E \geq \Big\lbrace \begin{array}{cc}
                       4\, \; \; & \mathrm{if} \; \; g=2 \\
                       2\, \; & \mathrm{if} \; \; g \geq3\\
                   \end{array}}$
     \end{itemize}
 \end{cor}

  \emph{Proof}.  In the situation of (a), $E$ must satisfy $G_d$, and the result was proved in \cite[4.6(c)]{SUV2003}. The remaining claims follow from \cref{GenGNN} with $d=5$ and $k=r(E)$, by noticing that assumption (d) is non-vacuous only when $\,\ell-e = 4\,$. $\blacksquare$\\

 Finally, from the proof of \cref{GenGNN} it follows that assumption (c) implies that $\,k \leq \ell-e-g+2\,$ (see also the discussion in \cite[page 14]{JU}) and that the Ext modules in assumption (d) automatically vanish when $g=2\,$ and $\,k = \ell-e-g+2$. Moreover, if $\,k \geq \ell-e-g+1\,$ we only need to assume $R$ to be Gorenstein locally in codimension $\ell-e$, since the Artin-Nagata property $AN_{\ell-e-k}$ is automatically satisfied when $\,\ell-e- k < g$. In this case, the depth condition in (c) can be simplified to assumption (c'). These observations prove the following two corollaries for modules with large reduction numbers, which recover \cite[6.5]{GNN} and \cite[6.4]{GNN} in the case when $E$ is an ideal.

 \begin{cor} 
   Let $R$ be a local Cohen-Macaulay ring with infinite residue field, and assume that $R$ is Gorenstein locally in codimension $\,\ell-e$. Let $E$ be a finite, torsion-free, orientable $R$-module, $\mathrm{rank}E=e>0$, $\ell(E)=\ell$, $\,\ell-e+1 \geq 2$. Let $g$ be the height of a generic Bourbaki ideal of $E$. Then, $\mathcal{R}(E)$ is Cohen-Macaulay if the following conditions hold:
     \begin{itemize}
       \item[$($a$)$] $E$ satisfies $G_{\,\ell-e+1}$. 
       \item[$($b$)$] $r(E) \leq \ell-e-g+1$. 
       \item[$($c$)$] $\,\mathrm{depth}_{\,}E^j \geq d-g-j+1\,$ for $\,1 \leq j \leq \ell-e-g+1$.
       \item[$($d$)$] If $g = 2$, $\,\mathrm{Ext}_{R_{\mathfrak{p}}}^{\,j+1}(E_{\mathfrak{p}}^j , R_{\mathfrak{p}}) =0\,$ for $\,1 \leq j \leq \ell-e-3, \,$ for all $\,\mathfrak{p} \in \mathrm{Spec}(R)\,$ with $\, \mathrm{dim} R_{\mathfrak{p}} = \ell-e\, $ such that $E_p$ is not free.
     \end{itemize}
  \end{cor}

 \begin{cor} \label{Lin3.4} \hypertarget{Lin3.4}{} 
   Let $R$ be a local Cohen-Macaulay ring with infinite residue field, and assume that $R$ is Gorenstein locally in codimension $\,\ell-e$. Let $E$ be a finite, torsion-free, orientable $R$-module, $\mathrm{rank}E=e>0$, $\ell(E)=\ell$, $\,\ell-e+1 \geq 2$. Let $g$ be the height of a generic Bourbaki ideal of $E$. Then, $\mathcal{R}(E)$ is Cohen-Macaulay if the following conditions hold:
    \begin{itemize}
       \item[$($a$)$] $E$ satisfies $G_{\,\ell-e+1}$. 
       \item[$($b$)$] $r(E) \leq \ell-e-g+2$. 
       \item[$($c$)$] $\mathrm{depth}_{\,}E^j \geq d-g-j+2\,$ for $\,1 \leq j \leq \ell-e-g+2$.
    \end{itemize}
 \end{cor}

  Notice that \cref{Lin3.4} recovers Lin's result \cite[3.4]{Lin}. In her proof, a generic Bourbaki ideal $I$ of $E$ is shown to satisfy the sliding depth property locally in codimension $\,\ell-e$. Recall that an $n$-generated ideal $I$ in a Cohen-Macaulay local ring $R$ satisfies \emph{sliding depth} if its Koszul homologies satisfy $\, \mathrm{depth}\,H_j \geq d-n+j\,$ for all $j \geq 1$ \cite{HVV}. Moreover, ideals with $G_{\infty}$ and sliding depth are of linear type by \cite[6.1]{HSVII}. In the proof of \cref{Lin3.4}, we only require $I$ to be of linear type locally in codimension $\,\ell-e$ instead. The improvement ultimately depends on being able to replace the depth condition of \cref{LinTypeDepth} with the weaker vanishing condition on the Ext modules in \cref{LinTypeExt}. \\

   
  We remark that the assumptions of \cref{GenGNN} are usually not satisfied for ideal modules. Nevertheless, by modifying the proof of \cref{GenGNN} we are able to provide a sufficient condition for the Rees algebra of an ideal module to be Cohen-Macaulay. This is discussed in the following theorem, which recovers Lin's \cite[4.3]{Lin}, with a simplified proof.
  
\begin{thm} \label{thmIdealModules} \hypertarget{thmIdealModules}{}
   Let $R$ be a local Cohen-Macaulay ring, and let $E$ be an ideal module with $\mathrm{rank}_{\,}E=e$ and $\ell(E)=\ell$. Assume that the following conditions hold.
    \begin{itemize}
        \item[$($a$)$] $r(E) \leq k$, where $k$ is an integer such that $1 \leq k \leq \ell -e$.
        \item[$($b$)$] $E$ is free locally in codimension $\,\ell-e- \mathrm{min} \{2, k\},\,$ and satisfies $G_{\,\ell-e+1}$.
        \item[$($c$)$]$\mathrm{depth}_{\,}E^j \geq d-\ell+e +k-j$ for $1 \leq j \leq k$.
     \end{itemize}
   Then, $\mathcal{R}(E)$ is Cohen-Macaulay.
 \end{thm}
  
  \emph{Proof}. Since $E$ is an ideal module, then $E$ admits a generic Bourbaki ideal $I$ with $\mathrm{ht}_{\,}I = g \geq 2$ and $r(I) \leq r(E) \leq k$, which satisfies $G_{\,\ell-e+1}$, that is, $G_{\,\ell(I)}$. Moreover, since $E$ is free locally in codimension $\,\ell-e- \mathrm{min} \{2, k\},\,$ by \cref{idealmodAN} $I$ satisfies $AN_{\,\ell-e-\mathrm{min} \{2, k\}}$. We next prove that 
   \begin{displaymath}
     \mathrm{depth}_{\,}I^j \geq d-\ell(I)+k-j+1= d-\ell +e +k-j \quad \mathrm{for} \quad 1 \leq j \leq k.
   \end{displaymath}
  Hence, by \cref{JUGNN} it would follow that $\, \mathcal{R}(I)$ is Cohen-Macaulay, whence the proof would be complete thanks to \cref{MainBourbaki}.
 
  The depth condition above is clearly satisfied if $\,e=1, \,$ so assume that $e\geq 2$. By \cref{DepthPasses} (with $s=d-\ell+e+k\,$) and \cref{SmallCodim} (with $s=\ell-e\,$), it suffices to show that the complexes $(\mathbb{C}'_j)_{\mathfrak{q}}$ are exact for all $\,\mathfrak{q} \in \mathrm{Spec}(R')\,$ with $\,\mathrm{dim}_{\,}R'_{\mathfrak{q}} \leq \ell-e\,$ and all $\, 0 \leq j \leq k$. 
  
  For each such $\mathfrak{q}$, let $\mathfrak{p}=\mathfrak{q} \cap R$. Then, $E_{\mathfrak{p}}$ is an ideal module, which is free locally in codimension $\,\ell-e- \mathrm{min} \{2, k\},\,$ and satisfies $G_{\,\ell-e+1}$, i.e. $G_{\infty}$. Hence, $\,E'_\mathfrak{p}/F'_\mathfrak{p}\,$ is isomorphic to an $R'_\mathfrak{p}$ ideal which satisfies $G_{\infty}$ and $AN_{\,\ell-e- \mathrm{min} \{2, k\}}$. Then by \cite[1.9]{Ulrich} its localization $\, (E'/F')_{\mathfrak{q}}\,$ satisfies $G_{\infty}$ and $AN_{\,\ell-e- \mathrm{min} \{2, k\}}$, whence also $AN_{\, \mathrm{dim}R'_{\mathfrak{q}} -2}\,$, since $\, \mathrm{dim}R'_{\mathfrak{q}} \leq \ell-e$. By \cite[1.9]{Ulrich} it then follows that $\,(E'/F')_{\mathfrak{q}}\,$ satisfies sliding depth. Hence, it is of linear type by \cite[6.1]{HSVII}, so that the complexes $\,(\mathbb{C}'_j)_{\mathfrak{q}}\,$ are exact for all $\, 0 \leq j \leq k$, thanks to \cref{SUV3.11}. $\blacksquare$\\

  Our last goal is to prove the following result, which is a module version of \cref{LinTypeDepth}. 
 
 \begin{thm} \label{herLinType} \hypertarget{herLinType}{}
    Let $R$ be a Gorenstein local ring of dimension $d$. Let $E$ be a finite, torsion-free, orientable $R$-module with $\mathrm{rank}_{\,}E =e >0$ and $\, \ell(E)=\ell$. Assume that $E$ satisfies $\,G_{\ell-e+2\,}$ and that $\, \mathrm{depth}_{\,}E^j \geq d-j\,$ for $1 \leq j \leq \ell-e-1$. Then, $E$ is of linear type and $\, \mathcal{R}(E) \,$ is Cohen-Macaulay. 
  \end{thm}
 
 This improves a result of Lin \cite[3.1]{Lin}, where the module was assumed to satisfy $G_{\infty}$ and the condition on the depths of powers $E^j$ was imposed on a larger range of exponents. In order to prove the theorem we will make use of the following preliminary lemma, which was inspired by the proof of \cite[3.1]{Lin}. 
 
 \begin{lemma} \label{IterDepthPasses} \hypertarget{IterDepthPasses}{}
   Let $R$ be a Cohen-Macaulay local ring of dimension $d$, and let $\,E= R a_1+ \ldots+ R a_n\,$ be a finite $R$-module. Denote
     $$\,\widetilde{R}\coloneq R[Z_1, \ldots, Z_n], \quad \widetilde{E} \coloneq E \otimes_R \widetilde{R} \quad \mathrm{and} \quad x \coloneq \sum_{i=1}^n Z_ia_i \in \widetilde{E}.$$
   For a positive integer $s$, assume that $\, \mathrm{depth}_{\,} E^j \geq d-j\,$ for $1 \leq j \leq s$. Then the following statements hold.
     \begin{itemize}
         \item[$($a$)$] For all $\, \mathfrak{q} \in \mathrm{Spec}\widetilde{R}\,$ and for all $\, 1 \leq j \leq s\,$, 
             $ \, \displaystyle{\mathrm{depth} \Big(\, \frac{\widetilde{E}_{\mathfrak{q}}^{\,j}}{x \widetilde{E}_{\mathfrak{q}}^{\,j-1}} \,\Big) \geq \mathrm{dim}_{\,} \widetilde{R}_{\mathfrak{q}} -j}$.
         \item[$($b$)$] Assume furthermore that $s \leq d-1$ and that $\, \displaystyle{(\widetilde{E}/\widetilde{R}x)_{\mathfrak{q}}}\,$ is of linear type for all $\, \mathfrak{q} \in \mathrm{Spec}\widetilde{R}\,$ such that $\, \mathrm{dim}_{\,}R_{\mathfrak{q} \cap R} \leq d-1. \,$ Then, for all $\,1 \leq j \leq s$
             $$\, \frac{\widetilde{E}^{\,j}}{x \widetilde{E}^{\,j-1}} \cong \Big(\, \frac{\widetilde{E}}{\widetilde{R}x} \,\Big)^j.$$
       \end{itemize}
 \end{lemma}
   
   \emph{Proof}. For $\, \mathfrak{q} \in \mathrm{Spec}(\widetilde{R})$, let $\mathfrak{p}= \mathfrak{q} \cap R$. Since $R_{\mathfrak{p}} \to \widetilde{R}_{\mathfrak{q}}$ is a flat local map, for all $1 \leq j \leq s$ we have that
       $$ \,  \mathrm{dim}_{\,}\widetilde{R}_{\mathfrak{q}} -\mathrm{depth}\,\widetilde{E}_{\mathfrak{q}}^j = \mathrm{dim}_{\,}R_{\mathfrak{p}} -\mathrm{depth}_{\,}E_{\mathfrak{p}}^{\,j} \leq d - \mathrm{depth}_{\,}E^j \leq j $$
   (see \cite[Theorem 15.1, Theorem 23.3 and Exercise 17.5]{Mats}). Moreover, by \cref{SUV3.6} $x$ is a regular element on $\,\mathcal{R}(\widetilde{E}),\,$ hence on $\,\mathcal{R}(\widetilde{E}_{\mathfrak{q}})$. Therefore, for all $j \geq 1$ we have $\,x\widetilde{E}_{\mathfrak{q}}^{j-1} \cong \widetilde{E}_{\mathfrak{q}}^{j-1}$. Hence, the conclusion in (a) follows from the Depth Lemma applied to the short exact sequence
      $$\, 0 \to x\widetilde{E}_{\mathfrak{q}}^{j-1} \longrightarrow \widetilde{E}_{\mathfrak{q}}^j \longrightarrow \frac{\widetilde{E}_{\mathfrak{q}}^{\,j}}{x \widetilde{E}_{\mathfrak{q}}^{\,j-1}} \to 0.$$
    
   Now assume that $s \leq d-1$. In order to prove (b), by \cref{SUV3.11} it suffices to show that $\, \widetilde{E}_{\mathfrak{q}}^j / x \widetilde{E}_{\mathfrak{q}}^{j-1} \,$ is a torsion-free $\widetilde{R}$-module for all $1 \leq j \leq s$. Equivalently, we need to show that $\, \displaystyle{\mathrm{depth} \Big(\, \frac{\widetilde{E}_{\mathfrak{q}}^{\,j}}{x \widetilde{E}_{\mathfrak{q}}^{\,j-1}} \,\Big) \geq 1 \,}$ for all $\,\mathfrak{q} \in \mathrm{Spec}(\widetilde{R})\,$ with $\, \mathrm{dim}_{\,}\widetilde{R}_{\mathfrak{q}} \geq 1$. 
  
   Let $\mathfrak{p}= \mathfrak{q} \cap R$. If $\mathrm{dim}_{\,}R_{\mathfrak{p}} \leq d-1, \,$ then by assumption $\, \displaystyle{(\widetilde{E}/\widetilde{R}x)_{\mathfrak{q}}}\,$ is of linear type. Hence, by \cref{SUV3.11} it follows that $\,\widetilde{E}_{\mathfrak{q}}^j / x \widetilde{E}_{\mathfrak{q}}^{j-1}\,$ is torsion-free for all $1 \leq j \leq s$. If $\mathrm{dim}_{\,}R_{\mathfrak{p}} = d, \,$ by part (a) we have that
    $$ \,\mathrm{depth} \Big(\, \frac{\widetilde{E}_{\mathfrak{q}}^{\,j}}{x \widetilde{E}_{\mathfrak{q}}^{\,j-1}} \,\Big) \geq \mathrm{dim}_{\,} \widetilde{R}_{\mathfrak{q}} -j \geq d-j,$$
   and the latter is at least 1 since $d \geq s+1$. This completes the proof. $\blacksquare$ \\

  \emph{Proof of Theorem 3.11}. Without loss of generality, we may assume that $E$ is not free. Let $E=Ra_1 + \ldots + Ra_n$. Since $E$ is torsion-free, orientable and satisfies $G_{\,\ell-e+2}$, by \cref{existBourbaki}, $\,E'/F' \cong J$ and $\,E''/F'' \cong I$, where $I$ and $J$ are ideals of height at least 2, satisfying $G_{\,\ell-e+2}$, i.e. $G_{\ell(I) +1}$ by \cref{anspreadBourbaki}. 
  
  By \cref{MainBourbaki}, it suffices to show that $I$ is of linear type and $\,\mathcal{R}(I)\,$ is Cohen-Macaulay. Let $g$ denote the height of $I$. If $e=1$, then $R =R''$ and $E \cong I$, so by assumption $\, \mathrm{depth}_{\,}I^j \geq d-j \geq d-j-g+2\,$ for $\,1 \leq j \leq \ell-e-1= \ell(I)-2,\,$ hence for $\,1 \leq j \leq \ell(I)-g$. Then, \cref{LinTypeDepth} implies that $I$ is of linear type and $\,\mathcal{R}(I)\,$ is Cohen-Macaulay.
  
  So, assume that $e\geq 2$. We induct on $\,d \geq 2$. If $\,d=2, \,$ then $\,\ell(I)=d=2=g\,$ and $I$ is $G_{\infty}$. This implies that $I$ is a complete intersection, thus $I$ satisfies sliding depth. Hence, by \cite[6.1]{HSVII} it follows that $I$ is of linear type and $\,\mathcal{R}(I)\,$ is Cohen-Macaulay. 
  
  Assume now that $d>2$. Notice that by \cref{LinTypeDepth} it suffices to show that 
   $$\, \mathrm{depth}_{\,}I^j \geq d-j \geq d-j-g+2 \quad \mathrm{for} \quad 1 \leq j \leq \ell-e-1= \ell(I)-2.$$
  By \cref{DepthPasses} and \cref{SmallCodim}, this follows once we show that the complexes $\,(\mathbb{C}'_j)_{\mathfrak{q}}\,$ are exact for all $\,\mathfrak{q} \in \mathrm{Spec}(R')\,$ with $\,\mathrm{dim}_{\,}R'_{\mathfrak{q}} \leq \ell-e\,$ and all $\, 0 \leq j \leq \ell-e-1$. In turn, by \cref{SUV3.11} we only need to show that $\, (E'/F')_{\mathfrak{q}} \cong J_{\mathfrak{q}}\, $ is of linear type for any such prime $\mathfrak{q}$. 
  
  For any such $\mathfrak{q}$, let $\mathfrak{p}= \mathfrak{q} \cap R$. Then, notice that $E_{\mathfrak{p}}$ is a finite, torsion-free, orientable $R_{\mathfrak{p}}$-module of rank $e$ and analytic spread $\,\ell(E_{\mathfrak{p}}) \leq \ell$. Moreover, our assumptions imply that
    $$ \, \mathrm{depth}_{\,}E_{\mathfrak{p}}^j \geq \mathrm{dim}_{\,}R_{\mathfrak{p}} -j \quad \mathrm{for} \quad 1 \leq j \leq \ell-e-1, $$
  hence for $\,1 \leq j \leq \ell(E_{\mathfrak{p}})-e-1$. Since $\,\mathrm{dim}_{\,}R_{\mathfrak{p}} \leq \mathrm{dim}_{\,}R'_{\mathfrak{q}} \leq \ell-e-1 \leq d-1, \,$ by the induction hypothesis $\,E_{\mathfrak{p}}\,$ is of linear type and $\,\mathcal{R}(E_{\mathfrak{p}}) \,$ is Cohen-Macaulay. 
  
  Now, let $\widetilde{R_{\mathfrak{p}}} \coloneq R_{\mathfrak{p}}[Z_1, \ldots, Z_n]$, $\,\displaystyle{\widetilde{E_{\mathfrak{p}}} \coloneq E_{\mathfrak{p}} \otimes_{R_{\mathfrak{p}}} \widetilde{R_{\mathfrak{p}}}}\,$ and $\, \displaystyle{ x \coloneq \sum_{i=1}^n Z_i \frac{a_i}{1}} \in \widetilde{E_{\mathfrak{p}}} \,$ be as in \cref{IterDepthPasses}. Since $\,E_{\mathfrak{p}}\,$ is of linear type and $\,\mathcal{R}(E_{\mathfrak{p}}) \,$ is Cohen-Macaulay, by \cref{MainBourbaki}(c) it follows  that $\, \displaystyle{(\widetilde{E_{\mathfrak{p}}}/\widetilde{R_{\mathfrak{p}}}x)}\,$ is of linear type. Hence, by \cref{IterDepthPasses} we have that 
   \begin{displaymath}
     \frac{\widetilde{E_{\mathfrak{p}}}^{\,j}}{x \widetilde{E_{\mathfrak{p}}}^{\,j-1}} \cong \Big(\, \frac{\widetilde{E_{\mathfrak{p}}}}{\widetilde{R_{\mathfrak{p}}}x} \,\Big)^j
   \end{displaymath}
  and that $ \, \displaystyle{\mathrm{depth} (\,\widetilde{E_{\mathfrak{p}}}/\widetilde{R_{\mathfrak{p}}}x \,) \geq \mathrm{dim}_{\,} \widetilde{R_{\mathfrak{p}}} -j} \,$ for all $\,1 \leq j \leq \ell-e-1$. Therefore, the same conclusions hold after tensoring with the ring $\, S\coloneq \widetilde{R_{\mathfrak{p}}}_{\,\mathfrak{m} \widetilde{R_{\mathfrak{p}}}}, \,$ where $\, \mathfrak{m}$ is the maximal ideal of $R$. In particular, the $S$-module $\, \displaystyle{(SE_{\mathfrak{p}}/Sx)}\,$ satisfies the same assumptions as the $R_{\mathfrak{p}}$-module $E_{\mathfrak{p}}$. Hence, we can iterate \cref{IterDepthPasses} $\,e-1$ times, to obtain that a generic Bourbaki ideal of $E_{\mathfrak{p}}$, constructed as in \cref{tdefBourbaki} with respect to the generators $\, \displaystyle \frac{a_1}{1}, \ldots, \frac{a_n}{1}, \,$ satisfies the inequalities
     $$\mathrm{depth} (\, E_{\mathfrak{p}}'' / F_{\mathfrak{p}}''\,) ^j \geq d-j \quad \mathrm{for} \quad 1 \leq j \leq \ell-e-1. \, $$
  The same ideal also satisfies $G_{\,\ell-e+2}$, hence $G_{\,\ell(E_{\mathfrak{p}}) -e+2}$, since $E_{\mathfrak{p}}$ does. Thus, by \cref{LinTypeDepth} it is of linear type, with Cohen-Macaulay Rees algebra. This implies that $E_{\mathfrak{p}}$ is of linear type, whence finally $J_{\mathfrak{q}}$ is of linear type by \cref{MainBourbaki}(c). $\blacksquare$ 
  

\end{document}